\documentclass[11pt,a4paper]{amsart}
\numberwithin{equation}{section}
\usepackage{amsmath, amssymb, amsthm}
\usepackage{mathtools}
\usepackage{mathrsfs}
\usepackage{array}
\usepackage{graphicx}
\usepackage{tikz}
\usepackage{comment}
\usepackage[margin=2.1cm]{geometry}
\allowdisplaybreaks[1]
\usepackage{hyperref}
\hypersetup{
  hidelinks,
  pdftitle={Hamiltonian Legendrian Curve Flows in Sasakian 3-Space Forms: mKdV equations and pseudo-Hermitian Geometry},
  pdfauthor={Shimpei Kobayashi}
}

\theoremstyle{plain}
\newtheorem{theorem}{Theorem}[section]
\newtheorem{lemma}[theorem]{Lemma}
\newtheorem{proposition}[theorem]{Proposition}
\newtheorem{corollary}[theorem]{Corollary}
\newtheorem{maintheorem}{Theorem}

\theoremstyle{definition}

\newtheorem{example}[theorem]{Example}

\theoremstyle{remark}
\newtheorem{remark}[theorem]{Remark}

\newcommand{\TW}{\mathrm{TW}}

\DeclareMathOperator{\sn}{sn}
\DeclareMathOperator{\cn}{cn}
\DeclareMathOperator{\dn}{dn}

\title[Hamiltonian Legendrian Curve Flows in Sasakian 3-Space Forms]
{Hamiltonian Legendrian Curve Flows in Sasakian 3-Space Forms:\\
 mKdV equations and pseudo-Hermitian Geometry}
\author{Shimpei Kobayashi}
\address{Department of Mathematics, Hokkaido University, Sapporo, Japan}
\email{shimpei@math.sci.hokudai.ac.jp}
\date{November 26, 2026}

\begin{document}
\begin{abstract}
We study the Hamiltonian flow of horizontal length on closed Legendrian
curves in three-dimensional Sasakian space forms. In the Heisenberg group,
the contact volume induces a symplectic form for the classical planar
filament flow on curves of zero signed area. Throughout the Sasakian
family, the Tanaka--Webster curvature satisfies the focusing modified
Korteweg--de Vries equation up to a translation in arclength.

We characterize the periodic curvature profiles that give closed curves
moving by ambient isometries in the Heisenberg group, and construct two families
in the universal cover of the real special linear group of degree two.
In both spaces the projected curve must close, and its horizontal lift
then closes precisely when the signed area vanishes. In the Heisenberg
group the zero-area locus is a single real-analytic graph joining the
translating figure-eight to the small-amplitude limit. Every rational rotation number strictly
between zero and one occurs as the turn per curvature period of a closed
rotating curve. In the hyperbolic model we obtain a
real-analytic family of closed curves moved by boosts and, for each
such rational turn, closed rotating curves at all sufficiently large
curvature scales. In both models, for every rational turn sufficiently
close to a full turn from below, a lifted rotation sweeps out an
embedded torus, at sufficiently large curvature scales in the hyperbolic
model.
\end{abstract}

\maketitle

\section{Introduction}
\label{sec:intro}

The classical planar filament flow moves a curve with normal speed $k_s$,
the arclength derivative of its curvature. Its curvature evolves by the
focusing modified Korteweg--de Vries (mKdV) equation
\cite{GoldsteinPetrich1991,NSW1992,ChouQu2002}. Langer and
Perline~\cite{LPPlanar} showed that length generates this motion through
a Poisson structure on plane curves. On a two-dimensional space form
the corresponding curvature equation acquires a constant transport term;
cf.~\cite{Anco2006}. In pseudo-Hermitian geometry, Calini, Ivey, and
Musso~\cite{CIM2024} obtained mKdV as the Hamiltonian flow of horizontal
length on Legendrian curves in the standard $3$-sphere. They also
constructed the hierarchy and studied its stationary curves.

We study the Hamiltonian flow of horizontal length on three-dimensional
Sasakian space forms. The contact volume induces a weakly symplectic form
on the smooth free locus of the reduced space of closed Legendrian curves.
In the Heisenberg group $H^1$, this form projects to the planar symplectic form
on curves of zero signed area. The planar formula uses an arclength
primitive of the normal variation. For a Legendrian lift, half the Reeb component
of the lifted variation is precisely such a primitive
(Remark~\ref{rem:planar-dictionary}). In every member of the Sasakian
family the Hamiltonian equation gives Reeb speed $2k$ up to a constant
along the curve, and the Legendrian condition then gives horizontal normal
speed $k_s$.

The Reeb projection sends a Legendrian curve in each Sasakian model
to a curve on a surface of constant Gaussian curvature. We write $G_\lambda$ for the
simply connected Sasakian $3$-space form whose Reeb quotient
$\Sigma_{2\lambda}$ has Gaussian curvature $2\lambda$. The
positive, zero and negative values of $\lambda$ give, respectively, a
Berger sphere (round at $\lambda=2$), the Heisenberg group $H^1$, and
$\widetilde{SL}(2,\mathbb R)$. For a Legendrian curve $\gamma$ with
projection $\beta$, the Reeb projection preserves arclength and
identifies the Tanaka--Webster curvature $k$ of $\gamma$ with the
geodesic curvature of $\beta$.
By Theorem~\ref{thm:sasakian-family-mkdv}, the Hamiltonian flow of
horizontal length evolves $k$ by
\begin{equation}
\label{eq:intro-mkdv}
 k_t=k_{sss}+\tfrac32 k^2k_s+2\lambda k_s.
\end{equation}
Changing the origin of
arclength by $2\lambda t$ removes the transport term. In particular,
the equation is focusing also in the pseudo-Hermitian anti-de Sitter
model, where a defocusing equation was anticipated in \cite[Section~5]{CIM2024}.

Which periodic traveling waves of \eqref{eq:intro-mkdv} give rise to
closed Legendrian curves? We study this question in the Heisenberg
group $H^1$ and in the universal cover $\widetilde{SL}(2,\mathbb R)$.
For a relative equilibrium, the curve evolves by ambient isometries up to
reparametrization, and its traveling-wave curvature profile $K(s-vt)$
satisfies the constrained-elastica equation on the base
\cite{LangerSinger1984,Heller2014}. Conversely, the profile determines
a Killing field on the base whose motion carries the reconstructed curve,
up to a shift in arclength. Its lift generates the relative equilibrium
by the momentum-map criterion of Theorem~\ref{thm:relative-equilibria}.
This recovers the motion from the curvature, but closure requires two
further tests. The isometry relating the projected Frenet frames one
curvature period apart must have finite order. Once the projection closes,
its horizontal lift may still end at a different point of the initial
Reeb fiber. In both noncompact models the displacement equals $-2$ times
the signed area of the projection (Theorem~\ref{thm:unified-holonomy}).
For instance, a Euclidean circle returns with its tangent after one turn,
yet its horizontal lift has a nonzero vertical displacement on each turn.
Vanishing signed area is therefore a closing condition in its own right.
It is preserved by the entire hierarchy
(Corollary~\ref{cor:lambda-holonomy-preservation}).

In $H^1$ we can rescale a nonconstant periodic profile so that its
curvature oscillates between $r$ and $1$. Write $\nu$ for the rescaled
wave speed. On the rotating branch, $\Theta(r,\nu)$ is the turning angle over
one curvature period, and $\mathcal B(r,\nu)$ is proportional to the
signed area swept by the radius vector over that period. Successive
arcs fit together after finitely many periods when $\Theta/(2\pi)$
is rational. Their lift closes when $\mathcal B=0$.

\begin{maintheorem}[The Heisenberg group]
\label{thm:intro-heisenberg}
\mbox{}
\begin{enumerate}
\item The zero set of $\mathcal B$ in the admissible region is the graph of
a single real analytic function $\nu=\nu_B(r)$, $-1\le r\le1$.  Its
endpoint $r=-1$ is the figure-eight elastica and $r=1$ is the
small-amplitude limit.  Along the graph the turning angle $\Theta$ runs
from $0$ to $2\pi$.
\item For every reduced fraction $p/n\in(0,1)$ there is a closed Legendrian
curve in $H^1$, closing after $n$ curvature periods, whose tangent turns
through $2\pi p/n$ in each period, and which the Hamiltonian flow of
horizontal length moves by rotations about the Reeb axis.
\item The curvature of every nonconstant closed relative equilibrium in
$H^1$ oscillates between the only two real roots of its quartic first
integral.  Among those moved by lifts of translations, the figure-eight is
the only one up to scaling and rigid motion.
\item Let $\Theta^*$ be the turning angle at $r=0$ on the graph; then
$\Theta^*<4\pi(1-1/\sqrt3)$.  If $2\pi p/n\ge\Theta^*$, the curve in
\textup{(2)} can be chosen so that rotation about the $z$-axis sweeps out an embedded torus.  In particular this holds for
$2-2/\sqrt3\le p/n<1$.
\end{enumerate}
\end{maintheorem}

\begin{figure}[!htbp]
 \centering
 \includegraphics[width=0.495\textwidth]{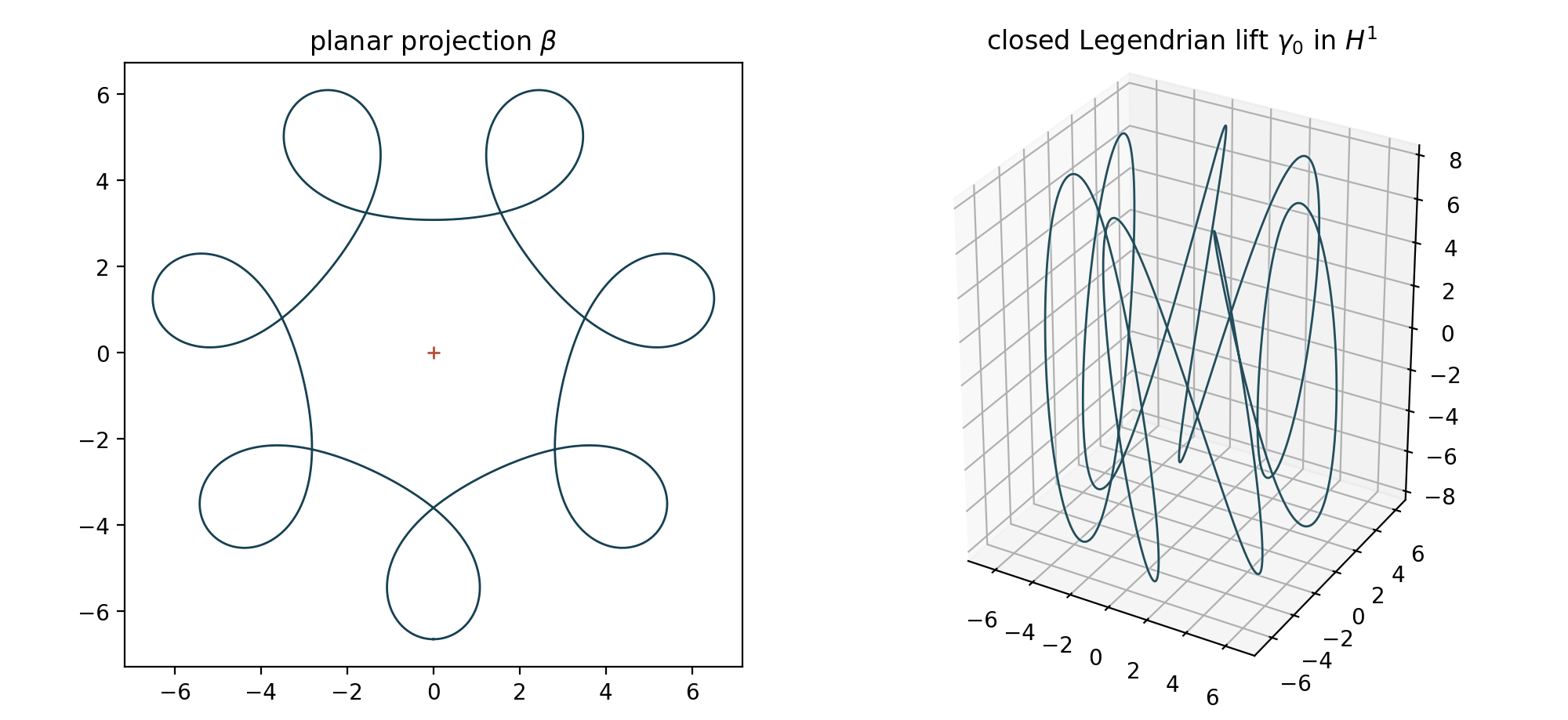}\hfill
 \includegraphics[width=0.495\textwidth]{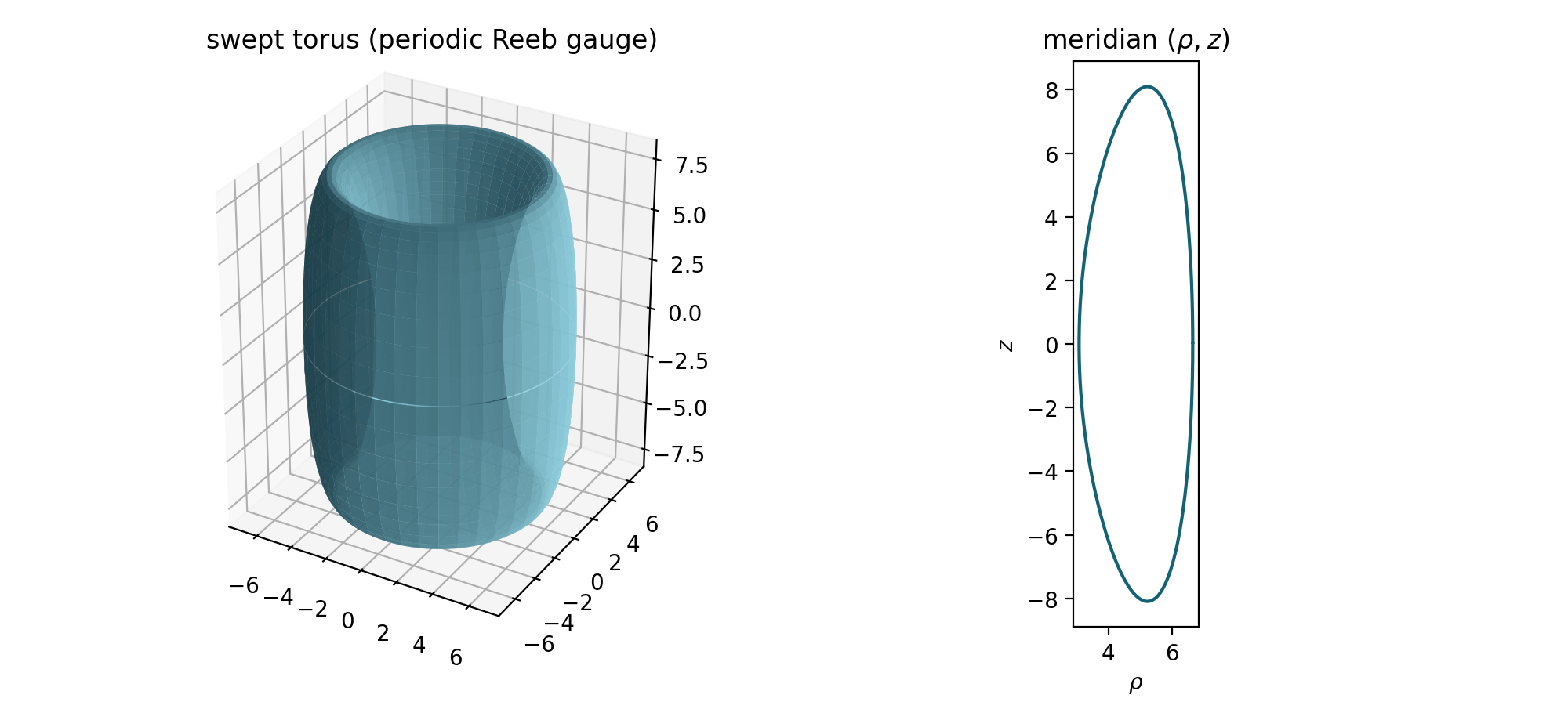}
 \caption{The rotating relative equilibrium of
 Corollary~\ref{cor:embedded-rotation-numbers} with \(p/n=6/7\), drawn with the
 numerical closing pair \(r=0.119457\ldots\), \(\nu=-0.163125\ldots\),
 \(c_\ast=0.050623\ldots\); the projected curve closes after seven curvature
 periods.  From left to right: the planar projection, the closed Legendrian
 lift, the embedded torus swept out by rotation about the $z$-axis, and its
 meridian, plotted in the unscaled group
 coordinates with equal coordinate units on all axes.}
 \label{fig:rotating-torus}
\end{figure}

Part~(1) is Lemma~\ref{lem:area-graph} with
Theorem~\ref{thm:closing-pairs}, part~(2) is
Theorem~\ref{thm:rotating-family}, part~(3) is
Proposition~\ref{prop:real-root-no-closure} with
Corollary~\ref{cor:translational-figure-eight}, and part~(4) is
Corollary~\ref{cor:embedded-rotation-numbers} with
Theorem~\ref{thm:rotating-torus}.
Figure~\ref{fig:flat-closing-map} shows how the two conditions meet.
The equation $\mathcal B=0$ determines a unique speed for each $r$.
As $r$ varies, the turning angle takes every value between $0$ and
$2\pi$, so prescribing a rational turn selects a closed rotation.
At the figure-eight endpoint the motion becomes a translation. There
$\mathcal B$ measures the displacement of the projected arc, and
reflection symmetry cancels the signed areas of the two lobes.
The translating lift and a rotating example are shown in
Figures~\ref{fig:flat-eight-motion} and~\ref{fig:rotating-torus}.

On the hyperbolic base, normalizing the ambient curvature to $-4$
leaves the curvature scale $b$ of the profile as a genuine parameter.
For the analogue of the figure-eight we determine the frequency for
which the projected frame returns after one curvature period. For
rotations the angular advance differs from the total curvature by an
area term. This correction vanishes on the zero-area locus, so the
same rational turning condition as in $H^1$ applies there
(Lemma~\ref{lem:split-rotating-data}).

\begin{maintheorem}[The special linear group]
\label{thm:intro-split}
Let $m_*\in(0,1)$ be the elliptic parameter with $2\mathbb E(m_*)=\mathbb K(m_*)$.
\begin{enumerate}
\item For every $m\in(m_*,1)$ there is a unique $\omega_m>0$ such that the
curvature $2\sqrt m\,\omega_m\cn(\omega_ms;m)$ reconstructs a closed
curve in $\mathbb H^2(-4)$; moreover $\omega_m>2$, and $m\mapsto\omega_m$ is real analytic.
For $m\in(0,m_*]$ no frequency gives a closed curve.  The
curve is the union of an arc and its image under a reflection, so its
signed area vanishes, and its horizontal lift is a closed Legendrian curve
in $\widetilde{SL}(2,\mathbb R)$ that the Hamiltonian flow moves by
hyperbolic boosts.  This motion has no nonzero time return.
\item For every reduced fraction $p/n\in(0,1)$ there is $b_0>0$ such that,
for every curvature scale $b>b_0$, there is a closed Legendrian curve in
$\widetilde{SL}(2,\mathbb R)$, closing after $n$ curvature periods, which
the flow moves by an elliptic one-parameter group.  If
$2\pi p/n>\Theta^*$, the lifted rotation fixing the fiber over the center pointwise sweeps out an embedded torus.
\end{enumerate}
\end{maintheorem}

Part~(1) is Theorem~\ref{thm:analytic-split-closure} with
Proposition~\ref{thm:cylinder-torus-dichotomy}, and part~(2) is
Theorems~\ref{thm:split-rotating} and~\ref{thm:split-torus}.
In part~(1), closure means that the boost displacement accumulated
\emph{along the curve} vanishes. The closed curve itself still moves
by a nontrivial boost \emph{in time} (Figure~\ref{fig:split-real-geometry}).
For part~(2), writing $\varepsilon=4/b^2$ gives a signed-area function
$\mathcal B_\varepsilon$ with $\mathcal B_0=\mathcal B$. The Heisenberg
zero-area graph persists on compact interior intervals for sufficiently
small $\varepsilon$. At larger $\varepsilon$ the numerical examples in
Figures~\ref{fig:split-parameter-map} and~\ref{fig:split-rotating-torus}
include closed rotating profiles with four real roots of the curvature
quartic, a type excluded in $H^1$ by Theorem~\ref{thm:intro-heisenberg}.

Closing in arclength and returning in time have different consequences
for the swept surface. Translations and boosts carry a closed curve
through a family with no time return. A lifted rotation fixing the
fiber over its center pointwise does return. To prove embeddedness one
must also compare the heights of the two points of the meridian at each
interior radius. This comparison gives the embedded tori in the two
theorems. Section~\ref{sec:surface} treats the local geometry of the
surface swept by the representative with Reeb speed $2k$. Where
$k\ne0$, the evolving curves are its characteristic curves and, with
the horizontal normal chosen there, its pseudo-Hermitian mean curvature
is $k$. Thus the same mKdV equation also governs this surface invariant.

\medskip\noindent\emph{Related work.}
Musso and Nicolodi~\cite{MussoNicolodi2026} construct the focusing
hierarchy by local motions of Legendrian curves in $H^1$. Closed
stationary solutions in the three pseudo-Hermitian models also feature
in the program of Calini, Ivey, Musso, and P\'ampano announced
in~\cite{MussoUFRJ}. For the projected curves, L.~Langer~\cite{Langer2025}
constructs planar zero-area elastica with prescribed rotation index,
and Heller~\cite{Heller2014} gives elliptic-function closing conditions
for constrained elastica in space forms. The rational datum in
Theorem~\ref{thm:intro-heisenberg} records the turn in a single curvature
period and hence the number of periods needed for the lift to close.
In anti-de Sitter geometry, the closed stationary curves of Musso and
P\'ampano~\cite{MussoPampano2024} are null for a Lorentzian metric.
The curves in Theorem~\ref{thm:intro-split} are horizontal for the
positive-definite Levi metric. Section~\ref{sec:cim-comparison}
compares the Hamiltonian motion with the Legendrian curve-shortening
flow of Pan and Sun~\cite{PanSun2021} and discusses the branches of
the hyperbolic closing problem that remain open.

\medskip\noindent\emph{Notation.}
We write $\sn(u;m)$, $\cn(u;m)$, and $\dn(u;m)$ for the Jacobi elliptic
functions with parameter $m$, and $\mathbb K(m)$, $\mathbb E(m)$, and
$\boldsymbol\Pi(\ell;m)$ for the complete elliptic integrals of the
first, second, and third kinds, respectively. The letter $K$ denotes a curvature profile and $w$ its
normalized form.  A Legendrian variation is written $\tau X+\zeta JX+\sigma\xi$, with tangential coefficient $\tau$, normal
coefficient $\zeta$ and Reeb coefficient $\sigma$; from
Section~\ref{sec:traveling-waves} on, $b$ is the curvature scale and $v$
the wave speed.  We write $\Omega_\lambda$ for the symplectic form and
$\mathcal{A}$ for the admissible parameter region.  The closing
functions are $\Theta(r,\nu)$ and $\mathcal B(r,\nu)$, the normalized
constants are $c_\ast(r,\nu)$ and $e_\ast(r,\nu)$, and $\mathscr D$ is a
signed filling chain.

\section{Sasakian \texorpdfstring{$3$}{3}-space forms and the Reeb fibration}
\label{sec:sasakian-family}

We set up the uniform model, the Legendrian Hamiltonian flow, and the
Reeb-holonomy closing law.  For the pseudo-Hermitian $3$-sphere, the
symplectic form, the Hamiltonian length flow, the mKdV hierarchy, and the
closing condition are due to Calini--Ivey--Musso~\cite{CIM2024}.  Their
constructions use only $\eta\wedge d\eta$ and the horizontal metric, and we
record them uniformly in the Sasaki $3$-space, with proofs only where the general case
differs.

\subsection{Models, curvature, and the Reeb quotient}
\label{sec:models-fibration}

For $\lambda\in\mathbb R$, let $G_\lambda$ be the simply connected Lie group
whose Lie algebra has a basis $(e_1,e_2,e_3)$ with
\begin{equation}
\label{eq:sasakian-brackets}
 [e_1,e_2]=-2e_3,\qquad
 [e_3,e_1]=-\lambda e_2,\qquad
 [e_2,e_3]=-\lambda e_1.
\end{equation}
Use the same symbols for the corresponding left-invariant vector fields,
and give $G_\lambda$ the metric $g$ for which this frame is orthonormal.
Put $\xi=e_3$ and $\eta=g(\xi,\cdot)$. By \eqref{eq:sasakian-brackets},
$d\eta(e_1,e_2)=2$ and $d\eta(\xi,\cdot)=0$. Thus $\eta$ is a contact
form with Reeb field $\xi$, and its horizontal plane field is
$H=\ker\eta=\operatorname{span}\{e_1,e_2\}$.

On each horizontal plane, let $J$ be rotation through $+\pi/2$ in
the orientation of $(e_1,e_2)$. This is the complex structure on $H$.
We extend it to $TG_\lambda$ by setting $J\xi=0$, so in our frame
\begin{equation}\label{eq:sasakian-J}
 Je_1=e_2,\qquad Je_2=-e_1,\qquad Je_3=0.
\end{equation}
In particular, if $X$ is the unit tangent to a Legendrian curve, then
$JX$ is its positively oriented horizontal unit normal.
Let $\nabla$ denote the Levi--Civita connection of $g$.
With the choice of $J$ in \eqref{eq:sasakian-J}, the Koszul formula gives
$\nabla_X\xi=JX$. The contact metric identities in this convention are
\[
 J^2=-I+\eta\otimes\xi,\qquad
 g(JX,JY)=g(X,Y)-\eta(X)\eta(Y),\qquad d\eta(X,Y)=2g(JX,Y).
\]
For background on Sasakian space forms and the Heisenberg model,
see~\cite{Blair2010,Boyer2009}.

\begin{proposition}[The Sasakian space forms]
\label{prop:sasakian-family}
The contact metric structure $(g,\eta,J,\xi)$ defined by
\eqref{eq:sasakian-brackets} and \eqref{eq:sasakian-J} is Sasakian,  i.e., 
$(\nabla_XJ)Y=-g(X,Y)\xi+\eta(Y)X$.
 Its sectional curvatures in the
distinguished planes are
\begin{equation}
\label{eq:sasakian-sectional-curvatures}
 K(e_1,e_2)=2\lambda-3,
 \qquad
 K(e_1,\xi)=K(e_2,\xi)=1.
\end{equation}
In particular, its constant $\phi$-sectional curvature is
$c_\phi=2\lambda-3$.  Moreover,
\[
 \mathfrak g_\lambda\cong
 \begin{cases}
   \mathfrak{su}(2),&\lambda>0,\\
   \mathfrak h_1,&\lambda=0,\\
   \mathfrak{sl}(2,\mathbb R),&\lambda<0.
 \end{cases}
\]
For $\lambda=2,0,-2$ one obtains respectively the round $3$-sphere $SU(2)$, the
standard Heisenberg group $H^1$, and a standard positive-definite left-invariant
Sasakian metric on the universal cover of the special linear group $\widetilde{SL}(2,\mathbb R)$.
\end{proposition}

This is standard; see \cite{Blair2010}.  In the frame
\eqref{eq:sasakian-brackets} the nonzero Levi--Civita derivatives are
$\nabla_{e_1}e_2=-\nabla_{e_2}e_1=-e_3$, $\nabla_{e_1}e_3=e_2$,
$\nabla_{e_2}e_3=-e_1$, $\nabla_{e_3}e_1=(1-\lambda)e_2$, and
$\nabla_{e_3}e_2=(\lambda-1)e_1$.

Let $\nabla^{\TW}$ denote the Tanaka--Webster connection
\cite{Tanaka1975,Webster1978} of this contact metric structure.  In the present sign convention it is related to the
Levi--Civita connection by
\begin{equation}
\label{eq:tw-lc-relation}
 \nabla^{\TW}_X Y
 =\nabla_XY-\eta(X)JY-\eta(Y)JX+g(JX,Y)\xi.
\end{equation}
For a Legendrian curve with unit horizontal tangent $X$, the
Tanaka--Webster curvature $k$ is defined by $\nabla^{\TW}_sX=kJX$.

\begin{remark}[The Heisenberg member in coordinates]
\label{rem:heisenberg-coordinates}
For $\lambda=0$ we use coordinates $(x,y,z)$ on $G_0=H^1\cong\mathbb R^3$
with
\begin{equation}
\label{eq:heisenberg-model}
 e_1=\partial_x+y\,\partial_z,\qquad
 e_2=\partial_y-x\,\partial_z,\qquad
 \xi=e_3=\partial_z,\qquad
 \eta=dz+x\,dy-y\,dx.
\end{equation}
Then $[e_1,e_2]=-2\xi$, $\eta(e_1)=\eta(e_2)=0$, and
$d\eta=2\,dx\wedge dy$, so that $d\eta(e_1,e_2)=2$.  The Reeb projection is
$(x,y,z)\mapsto(x,y)$, and a Legendrian curve $(x(s),y(s),z(s))$ satisfies
$z_s=-(xy_s-yx_s)$.  Thus the constant in $d\eta(e_1,e_2)=2$ agrees with the
normalization $A=2$ used for the Heisenberg group in
Chiu--Huang--Lai~\cite{CHL2017}, and the Tanaka--Webster curvature $k$ of a
Legendrian curve is the quantity called the $p$-curvature there.
\end{remark}

\begin{proposition}[The Reeb quotient and horizontal geometry]
\label{prop:reeb-quotient}
The Reeb action gives a Riemannian submersion
\[
 \pi_\lambda:G_\lambda\longrightarrow\Sigma_{2\lambda},
 \qquad d\pi_\lambda|_H:H\longrightarrow T\Sigma_{2\lambda}
 \ \text{is an isometry}.
\]
The base is the simply connected oriented surface of Gaussian curvature
$2\lambda$.  The compact model is a circle bundle; the simply connected
nilpotent and split models are $\mathbb R$-bundles.
A Legendrian curve is a horizontal lift of its projection, and its
Tanaka--Webster curvature is the signed geodesic curvature in the base.
\end{proposition}

\begin{proof}
The unit Reeb field is Killing and its flow preserves $H$ and its metric,
so the horizontal metric descends to the orbit space.  The standard
homogeneous quotients are the sphere, plane, and hyperbolic plane.
For horizontal $X,Y$, O'Neill's tensor is $A_XY=\tfrac12[X,Y]^{\rm vert}$
\cite{ONeill1966}, so $A_{e_1}e_2=-\xi$ by \eqref{eq:sasakian-brackets}.
O'Neill's curvature formula \cite[Corollary~1]{ONeill1966} and
\eqref{eq:sasakian-sectional-curvatures} give
\[
 K_\Sigma=K_{G_\lambda}(e_1,e_2)+3|A_{e_1}e_2|^2
 =(2\lambda-3)+3=2\lambda.
\]
The horizontal part of the Tanaka--Webster connection projects to the
base Levi--Civita connection.  Thus
$\nabla_s^{\TW}X=kJX$ projects to $\nabla_sT=kN$.
\end{proof}

\begin{remark}[D-homothety and the Berger spheres]
\label{rem:d-homothety}
For $\lambda>0$, $\lambda\ne2$, the metric $g$ is a Berger metric on
$S^3$; only $\lambda=2$ is round.  The Legendrian theory below, however,
uses only the horizontal metric $g_H$ and the contact form $\eta$: the
Tanaka--Webster curvature, the horizontal length \eqref{eq:horizontal-length},
the form \eqref{eq:sasakian-symplectic-form}, and the holonomy of
Theorem~\ref{thm:unified-holonomy}.  For $c>0$ the $D$-homothety
\[
 (g_H,\eta)\longmapsto(cg_H,c\eta),\qquad
 g=g_H+\eta\otimes\eta\longmapsto cg_H+c^2\eta\otimes\eta,
\]
preserves the normalization $d\eta=2\pi_\lambda^*dA$ and maps
$G_\lambda$ to $G_{\lambda/c}$.  It rescales $s\mapsto\sqrt c\,s$,
$k\mapsto k/\sqrt c$, $\mathcal L\mapsto\sqrt c\,\mathcal L$, and
$\Omega_\lambda\mapsto c^2\Omega_\lambda$, so the Hamiltonian field of
$\mathcal L$ is multiplied by $c^{-3/2}$ and \eqref{eq:sasakian-mkdv} is
preserved with $t\mapsto c^{3/2}t$.  Hence, for Legendrian curves, all
members with $\lambda>0$, including the Berger spheres, are equivalent to
the pseudo-Hermitian $3$-sphere of Calini--Ivey--Musso~\cite{CIM2024}, all
members with $\lambda<0$ are equivalent to one another, and $\lambda=0$ is
invariant.  The normalizations $\lambda=\pm2$ lose no generality, and for
$\lambda>0$ we refer to~\cite{CIM2024}.  Likewise, the universal cover
$A^3$ of anti-de Sitter space with its pseudo-Hermitian structure of
negative Webster curvature, listed in \cite[Section~1]{CIM2024}, is the case
$\lambda<0$: its Levi form is positive definite, so only the sign of
$\lambda$ distinguishes it.  The Lorentzian anti-de Sitter metric
$L_\theta-\theta^2$ and the Riemannian metric $L_\theta+\theta^2$ of
$G_{-2}$ share $(\theta,J)$ and the Levi form $L_\theta$, and differ only in
the Reeb direction, which Legendrian curves do not see.  With the
normalization of \cite{CIM2024}, in which the $3$-sphere has Webster
curvature $1$ and is our $G_2$, the space $A^3$ of Webster curvature $-1$ is
$G_{-2}$.
\end{remark}

\subsection{The reduced curve space and the length flow}
\label{sec:curve-space}

Let $\mathcal P$ be the space of Legendrian immersions
$\gamma\colon S^1=\mathbb R/\mathbb Z\to G_\lambda$, and let $\mathcal Q$
be its quotient by orientation-preserving reparametrizations and by the
Reeb action; manifold and weak-symplectic statements are made on its
smooth free locus.  When a closed curve is parametrized by horizontal
arclength $s$, its horizontal length is $\ell$.  For the orientation
of $(e_1,e_2,\xi)$, the Riemannian volume form of $g$ is
$\mu=\frac12\eta\wedge d\eta$, since $\mu(X,JX,\xi)=1$ for every
horizontal unit vector $X$.  Set
\begin{equation}
\label{eq:sasakian-symplectic-form}
 \Omega_\lambda(V,W)
 =-\frac14\int_{S^1}\mu(\gamma_x,V,W)\,dx .
\end{equation}
This is the form of \cite[(3.1)]{CIM2024} up to normalization.  For
$V=\tau X+\zeta JX+\sigma\xi$, Cartan's formula gives
$\partial_t|_{t=0}\eta(\gamma_s)=\partial_s(\eta(V))+d\eta(V,X)=\sigma_s-2\zeta$,
so the linearized Legendrian condition is
\begin{equation}
\label{eq:linearized-legendrian}
 \sigma_s=2\zeta,
\end{equation}
with no arclength condition imposed.  Tangential fields and constant
multiples of $\xi$ are the two quotient directions.  The horizontal length
\begin{equation}
\label{eq:horizontal-length}
 \mathcal L(\gamma)=\int_{S^1}\bigl|\gamma_x\bigr|\,dx
\end{equation}
is the pseudo-Hermitian length and descends to $\mathcal Q$.

\begin{lemma}[Legendrian variations]
\label{lem:legendrian-variations}
For a Legendrian variation $\gamma_t=\tau X+\zeta JX+\sigma\xi$ at fixed
parameter $x$, the arclength element $ds=|\gamma_x|\,dx$ satisfies
$\partial_t(ds)=(\tau_s-k\zeta)\,ds$. If $\partial_t(ds)=0$, then
\begin{equation}
\label{eq:sas-zc-a}
 \tau_s=k\zeta,\qquad
 \sigma_s=2\zeta,\qquad
 k_t=(\zeta_s+k\tau)_s+2\lambda\zeta.
\end{equation}
\end{lemma}

\begin{proof}
Equation~\eqref{eq:linearized-legendrian} gives $\sigma_s=2\zeta$.
For the projection $\beta$ of speed $q=|\gamma_x|$,
Proposition~\ref{prop:reeb-quotient} and the Frenet equations give
$q_t=\langle\nabla_x^\Sigma\beta_t,T\rangle
=q\langle\nabla_s^\Sigma(\tau T+\zeta N),T\rangle=q(\tau_s-k\zeta)$.
Thus $ds=q\,dx$ has the stated variation. When $\tau_s=k\zeta$, the arclength
curvature variation formula on $\Sigma_{2\lambda}$ gives \eqref{eq:sas-zc-a}.
\end{proof}

\begin{theorem}[Horizontal length generates focusing mKdV]
\label{thm:sasakian-family-mkdv}
On the smooth free locus of $\mathcal Q$, the two-form $\Omega_\lambda$
induced by $\mu$ is weakly symplectic. The Hamiltonian field of horizontal length has
normal speed $\zeta=k_s$ and Reeb speed $\sigma=2k$ modulo constants in $s$.
An arclength-preserving representative is
\begin{equation}
\label{eq:sasakian-curve-flow}
 \gamma_t=\frac12k^2X+k_sJX+2k\xi,
\end{equation}
and the Tanaka--Webster curvature satisfies
\begin{equation}
\label{eq:sasakian-mkdv}
 k_t=k_{sss}+\frac32k^2k_s+2\lambda k_s.
\end{equation}
Thus $\widehat k(s,t)=k(s-2\lambda t,t)$ solves focusing mKdV, recovering
\cite[Theorem~3.2, Proposition~3.3]{CIM2024} at $\lambda=2$.
\end{theorem}

\begin{proof}
As in \cite{CIM2024}, $\Omega_\lambda$ is the transgression of the
closed form $\mu$ \cite{Vizman2011,KM1997}. For Legendrian variations
$V,W$, write $\sigma=\eta(V)$ and $\widetilde\sigma=\eta(W)$ for their Reeb
components. By \eqref{eq:linearized-legendrian} and integration by parts,
\begin{equation}
\label{eq:omega-reeb}
 \Omega_\lambda(V,W)=\frac14\int_{S^1}\sigma\,\widetilde\sigma_s\,ds ,
\end{equation}
which vanishes for all $W$ only if $\sigma$ is constant.
For $V=\tau X+\zeta JX+\sigma\xi$, Lemma~\ref{lem:legendrian-variations} gives
$\delta_V(ds)=(\tau_s-k\zeta)\,ds$. Since $\tau$ is periodic,
$\int_0^\ell \tau_s\,ds=0$, so $d\mathcal L(V)=-\int k\zeta\,ds$.
Using \eqref{eq:linearized-legendrian} and integrating by parts gives
$d\mathcal L(V)=\frac12\int k_s\sigma\,ds=\Omega_\lambda(V,k_sJX+2k\xi)$.
The periodic choice $\tau=k^2/2$ gives $\tau_s=k\zeta$ for $\zeta=k_s$. Thus
\eqref{eq:sas-zc-a} yields \eqref{eq:sasakian-mkdv}.
\end{proof}

For $\lambda=2,0,-2$ ($c_\phi=1,-3,-7$) the transport term in
\eqref{eq:sasakian-mkdv} is $+4k_s$, $0$, and $-4k_s$; the first agrees with
\cite{CIM2024}.  In particular the case $\lambda<0$, i.e.\ the
pseudo-Hermitian $A^3$ of \cite{CIM2024}, also gives focusing mKdV;
\cite[Section~5]{CIM2024} expected the defocusing equation there.  The
reason is visible in \cite[(3.7)]{CIM2024}, in our notation $k_t=\zeta_{ss}+(k\tau)_s+4\zeta$: the ambient geometry enters
only through the zeroth-order term $4\zeta$, which becomes $2\lambda\zeta$, while
the cubic term comes from $(k\tau)_s$ with $\tau_s=k\zeta$.  Since the Levi form of
a three-dimensional pseudo-Hermitian structure is definite, a defocusing
Legendrian equation would require leaving CR geometry, for instance a
Lorentzian horizontal metric.

\begin{remark}[The planar mKdV flow]
\label{rem:planar-dictionary}
Consider oriented closed plane curves of fixed signed area, modulo
reparametrization. The normal component $\zeta_V$ of a tangent vector $V$
satisfies $\oint\zeta_V\,ds=0$, so it has a periodic primitive $\chi_V$ with
$(\chi_V)_s=\zeta_V$. The planar form and the first variation of length are
\[
 \Omega_{\mathrm{pl}}(V,W)=\oint \chi_V\zeta_W\,ds,\qquad
 d\mathcal L(V)=-\oint k\zeta_V\,ds=\oint \chi_V k_s\,ds .
\]
The additive constant in $\chi_V$ drops out because $\oint\zeta_W\,ds=0$.
Thus $d\mathcal L(V)=\Omega_{\mathrm{pl}}(V,X_{\mathcal L})$ gives
$\zeta_{X_{\mathcal L}}=k_s$, the normal speed of the planar filament
flow~\cite{LPPlanar}.

In the Heisenberg group a closed plane curve admits a closed Legendrian
lift precisely when its signed area vanishes
(Theorem~\ref{thm:unified-holonomy}). For a lifted variation $\widehat V$
with Reeb component $\sigma_V$, the Legendrian condition reads
$(\sigma_V)_s=2\zeta_V$. Taking $\chi_V=\sigma_V/2$ in the planar form gives
$\Omega_{\mathrm{pl}}(V,W)=\Omega_0(\widehat V,\widehat W)$ by
\eqref{eq:omega-reeb}. Adding a constant Reeb field to $\widehat V$
accounts for the freedom to choose the constant in $\chi_V$.
\end{remark}

\begin{corollary}[Projection of the Hamiltonian flow]
\label{thm:sasaki-projection}
The flow \eqref{eq:sasakian-curve-flow} projects to
$\beta_t=\frac12k^2T+k_sN$ on $\Sigma_{2\lambda}$, and its geodesic
curvature satisfies \eqref{eq:sasakian-mkdv}.
\end{corollary}

\begin{theorem}[Legendrian holonomy and closure]
\label{thm:unified-holonomy}
Let $\beta:S^1\to\Sigma_{2\lambda}$ be an oriented immersed closed curve and
$\mathscr D$ an integral two-chain with $\partial\mathscr D=\beta$, area
counted with multiplicity.  Its horizontal lift has Reeb displacement
\begin{equation}
\label{eq:unified-holonomy}
 \operatorname{Hol}(\beta)
 =\exp\!\left(-2\operatorname{Area}_{\rm sgn}(\mathscr D)\,\xi\right).
\end{equation}
For a real-line fiber, the lift closes precisely when the signed area
vanishes. For a circle fiber, let $P$ be the least positive return time
of the Reeb flow $\Phi_t^\xi$ at a point $x$ of the fiber, so that
$\Phi_P^\xi(x)=x$. The lift then closes precisely when
\begin{equation}
\label{eq:quantized-area}
 2\operatorname{Area}_{\rm sgn}(\mathscr D)\in P\mathbb Z.
\end{equation}
Moreover,
\begin{equation}
\label{eq:gb-holonomy}
 \int_\beta k\,ds+2\lambda\operatorname{Area}_{\rm sgn}(\mathscr D)
 =2\pi\operatorname{rot}_{\mathscr D}(\beta).
\end{equation}
\end{theorem}

\begin{proof}
The curvature of the connection is $2\,dA_{2\lambda}$ on the base.
By Stokes' theorem the Reeb displacement is
$-2\operatorname{Area}_{\rm sgn}(\mathscr D)$ modulo the fiber period
when the fiber is compact, and exactly this value otherwise.
Equation~\eqref{eq:gb-holonomy} follows from Gauss--Bonnet.
\end{proof}

For $\lambda>0$ the Reeb period is $P=4\pi/\lambda$ and the base sphere has
area $2\pi/\lambda$, so \eqref{eq:quantized-area} is independent of the
filling. By \eqref{eq:gb-holonomy}, the horizontal lift closes after
one traversal of $\beta$ precisely when
\[
 \operatorname{rot}_{\mathscr D}(\beta)
 -\frac{1}{2\pi}\int_\beta k\,ds\ \in 2\mathbb Z .
\]
In particular a closed lift has integral normalized total curvature,
as in \cite[Proposition~2.8]{CIM2024}. The parity matters: a great circle
traversed once has $k=0$ and bounds a hemisphere of area $\pi/\lambda$.
Its lift acquires half a Reeb period and closes on the second traversal.  In $H^1$ and
$\widetilde{SL}(2,\mathbb R)$ the condition is zero signed area.

\subsection{The shifted hierarchy}
\label{sec:hierarchy}

Let $\partial_s^{-1}$ denote a periodic primitive of a mean-zero function,
or the local primitive displayed by a differential-polynomial formula, and
let $\mathcal R=\partial_s^2+k^2+k_s\partial_s^{-1}k$ be the focusing mKdV
recursion operator~\cite{Olver1993}.  For $\lambda=2$ the hierarchy below
is \cite[Section~3.2]{CIM2024}.

\begin{lemma}[Hamiltonian normal coefficients]
\label{lem:hamiltonian-normal}
Let $\mathcal H$ be an invariant functional on $\mathcal Q$ with
$d\mathcal H(V)=\int_0^\ell G_{\mathcal H}\zeta\,ds$.  Its Hamiltonian field,
$d\mathcal H(V)=\Omega_\lambda(V,X_{\mathcal H})$, has Reeb coefficient
$-2G_{\mathcal H}$ and normal coefficient
$\zeta_{\mathcal H}=-\partial_sG_{\mathcal H}$, and
$d\mathcal H(X_{\mathcal K})=-\int_0^\ell G_{\mathcal H}\,\partial_sG_{\mathcal K}\,ds$.
\end{lemma}

\begin{proof}
For $W=X_{\mathcal H}$ insert $\widetilde\sigma=-2G_{\mathcal H}$ in \eqref{eq:omega-reeb} and use
\eqref{eq:linearized-legendrian}.
\end{proof}

For $\mathcal H=\int f\,ds$ with $f$ a differential polynomial in $k$, let
$h=\delta\mathcal H/\delta k$ with $ds$ held fixed and let $S_f$ be the
stress, $\partial_sS_f=hk_s$, with constant fixed by $f$.  Since
$\delta(ds)=(\tau_s-k\zeta)\,ds$ and
$\delta k=\zeta_{ss}+(k^2+2\lambda)\zeta+\tau k_s$, the normal gradient is
\begin{equation}
\label{eq:normal-gradient-stress}
 G_{\mathcal H}=(\partial_s^2+k^2+2\lambda)h-kS_f .
\end{equation}

\begin{theorem}[Shifted hierarchy on $G_\lambda$]
\label{thm:lambda-hierarchy}
Let $P_0=k$, $\zeta_j=\partial_sP_j$, $\zeta_{j+1}=\mathcal R\zeta_j$, $\mathcal
H_0=\mathcal L$, and for $j\ge1$ let $\mathcal H_j$ be the standard mKdV
integral with $\delta\mathcal H_j/\delta k=-P_{j-1}$.  The Hamiltonian field
of $\mathcal H_j$ on $(\mathcal Q,\Omega_\lambda)$ has normal coefficient
\begin{equation}
\label{eq:lambda-Hj-normal}
 \zeta_j+2\lambda\zeta_{j-1}\qquad(j\ge1),
\end{equation}
and curvature evolution $k_t=\mathcal R_\lambda(\zeta)$ with
$\mathcal R_\lambda=\mathcal R+2\lambda I$.  Consequently
\begin{equation}
\label{eq:shifted-hamiltonians}
 \mathcal K_0^{(\lambda)}=\mathcal H_0,\qquad
 \mathcal K_j^{(\lambda)}=\sum_{r=1}^j\binom{j-1}{r-1}(2\lambda)^{j-r}\mathcal H_r
\end{equation}
has normal coefficient $\mathcal R_\lambda^j(k_s)$ and induces
$k_{t_j}=\mathcal R_\lambda^{j+1}(k_s)$.  These Hamiltonians are in
involution, and their flows commute.
\end{theorem}

\begin{proof}
Lemma~\ref{lem:legendrian-variations} with $\tau=\partial_s^{-1}(k\zeta)$ gives
$k_t=\mathcal R_\lambda\zeta$.  The Lenard identity
$(\partial_s^2+k^2)h_j-kS_j=-P_j$ \cite{Olver1993} and
\eqref{eq:normal-gradient-stress} give
$G_{\mathcal H_j}=-(P_j+2\lambda P_{j-1})$, hence
\eqref{eq:lambda-Hj-normal} by Lemma~\ref{lem:hamiltonian-normal}.  The
binomial identity
$\binom{j-1}{r-1}+\binom{j-1}{r}=\binom jr$ turns
\eqref{eq:shifted-hamiltonians} into $\sum_r\binom jr(2\lambda)^{j-r}\zeta_r
=\mathcal R_\lambda^j(k_s)$.  Involution follows from
$\int P_r\partial_sP_q\,ds=0$ and Lemma~\ref{lem:hamiltonian-normal}.
\end{proof}

\begin{corollary}[The hierarchy preserves Reeb holonomy]
\label{cor:lambda-holonomy-preservation}
Every flow of the hierarchy preserves the signed area of the projected curve,
hence the closing condition \eqref{eq:quantized-area}.
\end{corollary}

\begin{proof}
The normal velocity is a combination of the total derivatives
$\zeta_r=\partial_sP_r$, so $\oint\zeta\,ds=0$.
\end{proof}

\section{Global reconstruction, relative equilibria, and swept topology}
\label{sec:global-reconstruction}

We characterize relative equilibria by momentum maps, derive their
traveling-wave ODE, and reconstruct the Killing generator on the Reeb
quotient.  The closing problem then consists of finite-order Frenet
monodromy in the base and the Reeb holonomy of
Theorem~\ref{thm:unified-holonomy}.  We solve it explicitly for the flat and
hyperbolic bases; for the spherical base see~\cite{CIM2024}.

\subsection{Symmetries, momentum maps, and relative equilibria}
\label{sec:relative-equilibria-criterion}
A relative equilibrium is a closed Legendrian curve whose Hamiltonian
motion agrees, in the reduced curve space, with an ambient Sasakian
isometry.  For the $3$-sphere these are the stationary curves of
\cite[Section~4]{CIM2024}, where the conserved momentum \cite[(4.6)]{CIM2024}
plays the role of the Killing data $(Y,\psi)$ of
Proposition~\ref{prop:space-form-noether-reconstruction} below, and closure
is characterized by the rationality of two periods
\cite[Proposition~4.3, Remark~4.4]{CIM2024}.  The detailed analysis there
concerns the phase-symmetric branch, which corresponds to $C=0$ below.

Let $Y$ be a Killing field of $\Sigma_{2\lambda}$.  A \emph{lift} of $Y$ is
an infinitesimal automorphism $\widehat Y$ of $(g,\eta,J,\xi)$ with
$d\pi_\lambda(\widehat Y)=Y$.  Since $\mathcal L_{\widehat Y}\eta=0$,
Cartan's formula gives $d(\eta(\widehat Y))=-\iota_{\widehat Y}d\eta$; as
$\iota_\xi d\eta=0$, the function $\eta(\widehat Y)$ is constant along the
Reeb orbits, and $d\eta=2\pi_\lambda^*dA_{2\lambda}$ gives
\begin{equation}
\label{eq:killing-potential}
 \eta(\widehat Y)=f_Y\circ\pi_\lambda,\qquad df_Y=-2\,\iota_Y dA_{2\lambda}.
\end{equation}
We call $f_Y$ the \emph{Killing potential} of $\widehat Y$.  Conversely,
$Y$ preserves $dA_{2\lambda}$ and $\Sigma_{2\lambda}$ is simply connected,
so a function $f_Y$ with \eqref{eq:killing-potential} exists and is unique
up to an additive constant; the field $Y^{h}+(f_Y\circ\pi_\lambda)\xi$,
with $Y^h$ the horizontal lift, preserves $\eta$ by
\eqref{eq:killing-potential}, and it preserves the horizontal metric and
$J|_H$ because these are pulled back from the base.  Adding a constant to
$f_Y$ adds a constant multiple of $\xi$ to $\widehat Y$.

\begin{example}
\label{ex:heisenberg-potentials}
In the coordinates \eqref{eq:heisenberg-model} of $H^1$, the rotation
$Y=-y\,\partial_x+x\,\partial_y$ has potential $f_Y=x^2+y^2$, and its lift
is the rotation about the $z$-axis.  The translation
$Y_w=w_1\partial_x+w_2\partial_y$ has potential
$f_{Y_w}(x,y)=2(w_2x-w_1y)=2\det\bigl((x,y),w\bigr)$, and its lift is a left
translation of $H^1$.  On $\mathbb H^2(2\lambda)$ and $S^2(2\lambda)$ the
Killing potentials are, up to constants and a normalization factor, the
restrictions of linear functions of the hyperboloid or sphere model.
\end{example}

\begin{proposition}[Momentum maps]
\label{prop:momentum-maps}
Let $\widehat Y$ be a lift of a Killing field $Y$ with potential $f_Y$.  For
$[\gamma]\in\mathcal Q$ with projection $\beta$ and an integral two-chain
$\mathscr D$ with $\partial \mathscr D=\beta$, put
\begin{equation}
\label{eq:momentum-map}
 \mathcal M_Y([\gamma])=\frac12\int_{\mathscr D}f_Y\,dA_{2\lambda}.
\end{equation}
For $\lambda\le0$ the integral is independent of the compact signed
filling.  For $\lambda>0$ normalize the lift by a constant Reeb field
(which acts trivially on $\mathcal Q$), and choose the mean-zero potential
$\int_{\Sigma_{2\lambda}}f_Y\,dA_{2\lambda}=0$; the integral is then also
filling-independent, since two fillings differ by a multiple of the sphere.
Other potential normalizations give local Hamiltonians differing by
constants on each fixed-area branch.  Then $\widehat Y$ restricted to $\gamma$ is
an admissible variation, and
\begin{equation}
\label{eq:momentum-hamiltonian}
 d\mathcal M_Y(V)=\Omega_\lambda(V,\widehat Y)
\end{equation}
for every admissible $V$.  Thus $\widehat Y$ is the Hamiltonian vector field
of $\mathcal M_Y$ on $\mathcal Q$, with normal gradient
$G_{\mathcal M_Y}=-\frac12f_Y\circ\beta$.  Moreover $\mathcal M_Y$ is
conserved by every flow generated by $\mathcal K_j^{(\lambda)}$.
\end{proposition}

\begin{proof}
The flow of $\widehat Y$ consists of Sasakian automorphisms, which carry
Legendrian curves to Legendrian curves, so $\widehat Y|_\gamma$ is
admissible, and its Reeb coefficient is $\eta(\widehat Y)=f_Y\circ\beta$ by
\eqref{eq:killing-potential}.  Let $V=\tau X+\zeta JX+\sigma\xi$.  Its projection moves
$\beta$ with normal velocity $\zeta$ along $N=JT$, and $N=JT$ is the left normal for the boundary orientation.  The
boundary-variation formula with signed multiplicities therefore gives the first variation of \eqref{eq:momentum-map} as
$d\mathcal M_Y(V)=-\frac12\int_0^\ell f_Y(\beta)\,\zeta\,ds$.  On the other hand,
by \eqref{eq:omega-reeb} and \eqref{eq:linearized-legendrian},
\[
 \Omega_\lambda(V,\widehat Y)
 =\frac14\int_0^\ell \sigma\,\bigl(f_Y\circ\beta\bigr)_s\,ds
 =-\frac14\int_0^\ell \sigma_s\,f_Y(\beta)\,ds
 =-\frac12\int_0^\ell \zeta\,f_Y(\beta)\,ds,
\]
which proves \eqref{eq:momentum-hamiltonian}.  Every
$\mathcal K_j^{(\lambda)}$ is an integral of a differential polynomial in
the curvature, or the length, and is therefore invariant under the flow of
$\widehat Y$; hence
$d\mathcal M_Y(X_{\mathcal K})=\Omega_\lambda(X_{\mathcal K},\widehat Y)
=-d\mathcal K(\widehat Y)=0$.
\end{proof}

\begin{theorem}[Variational characterization of relative equilibria]
\label{thm:relative-equilibria}
Let $\widehat Y$ be a lift of a Killing field $Y$ with potential $f_Y$, and
let $[\gamma]\in\mathcal Q$ have projection $\beta$ and curvature $k$.  The
following are equivalent:
\begin{enumerate}
\item\label{it:re-field} $X_{\mathcal L}([\gamma])=\widehat Y([\gamma])$ in
$T_{[\gamma]}\mathcal Q$;
\item\label{it:re-critical} $[\gamma]$ is a critical point of
$\mathcal L-\mathcal M_Y$ on $\mathcal Q$;
\item\label{it:re-local} $2k-f_Y\circ\beta$ is constant along $\beta$.
\end{enumerate}
In this case $t\mapsto\exp(t\widehat Y)\cdot[\gamma]$ is an integral curve
of $X_{\mathcal L}$ on $\mathcal Q$.  More precisely, if
$2k-f_Y\circ\beta=-c_0$, then the representative
\begin{equation}
\label{eq:relative-equilibrium-velocity}
 \gamma_t=\frac12k^2X+k_sJX+(2k+c_0)\xi
\end{equation}
of $X_{\mathcal L}$ differs from $\widehat Y|_\gamma$ by a tangential field,
so that $\exp(t\widehat Y)\circ\gamma$, reparametrized by horizontal
arclength, satisfies \eqref{eq:relative-equilibrium-velocity} with this value of $c_0$.
\end{theorem}

\begin{proof}
A tangent vector of $\mathcal Q$ is represented by an admissible field
modulo tangential fields and constant multiples of $\xi$.  Since the normal
coefficient is $\zeta=\frac12\sigma_s$ by \eqref{eq:linearized-legendrian}, the
class is determined by the Reeb coefficient $\sigma$ modulo constants.  The
Reeb coefficient of $X_{\mathcal L}$ is $2k$ by
Theorem~\ref{thm:sasakian-family-mkdv}, and that of $\widehat Y$ is
$f_Y\circ\beta$ by \eqref{eq:killing-potential}.  This proves
(\ref{it:re-field})$\Leftrightarrow$(\ref{it:re-local}) and the statement about
\eqref{eq:relative-equilibrium-velocity}.  By the first variation formula in
the proof of Theorem~\ref{thm:sasakian-family-mkdv} and
Proposition~\ref{prop:momentum-maps},
\[
 d(\mathcal L-\mathcal M_Y)(V)
 =-\int_0^\ell\Bigl(k-\tfrac12f_Y\circ\beta\Bigr)\zeta\,ds.
\]
The admissible normal coefficients are exactly the functions $\zeta$ with
$\oint\zeta\,ds=0$, because $\zeta=\frac12\sigma_s$ with $\sigma$ periodic.  Hence
(\ref{it:re-critical}) holds if and only if $k-\frac12f_Y\circ\beta$ is
constant, which is (\ref{it:re-local}).

Finally, let $\phi_t=\exp(t\widehat Y)$.  It preserves $\mathcal L$ and
$\Omega_\lambda$, hence $d\phi_t\,X_{\mathcal L}=X_{\mathcal L}\circ\phi_t$, and
\[
 \frac{d}{dt}\phi_t[\gamma]
 =\widehat Y(\phi_t[\gamma])
 =d\phi_t\,\widehat Y([\gamma])
 =d\phi_t\,X_{\mathcal L}([\gamma])
 =X_{\mathcal L}(\phi_t[\gamma]).
\qedhere
\]
\end{proof}

\begin{remark}[The role of the constant $c_0$]
\label{rem:relative-equilibrium-constant}
The constant in (\ref{it:re-local}) plays two roles.  Variationally, it is the
Lagrange multiplier of the fixed signed-area class, since $\mathcal Q$ only
admits area-preserving variations.  Geometrically, adding a constant multiple of $\xi$ to the velocity
gives the same tangent vector on $\mathcal Q$. For the chosen lift
$\widehat Y$, taking $c_0=f_Y\circ\beta-2k$ in
\eqref{eq:relative-equilibrium-velocity} makes the curve follow the orbit of
$\exp(t\widehat Y)$ up to reparametrization. The swept surface can depend
on this choice (Remark~\ref{rem:torus-reeb-component}).  Since isometries preserve
curvature, the curvature of every relative equilibrium is a traveling wave
in the sense of Proposition~\ref{prop:traveling-wave}.
\end{remark}

\subsection{General traveling waves on
\texorpdfstring{$G_\lambda$}{G-lambda}}
\label{sec:traveling-waves}

Write $Y=\tau T+\zeta N$ along $\beta$, where $T=\beta_s$ and $N=JT$. By
\eqref{eq:killing-potential}, $(f_Y\circ\beta)_s=-2dA(Y,T)=2\zeta$ since $dA(T,N)=1$.
At a relative equilibrium $f_Y\circ\beta=2k+c_0$ by Theorem~\ref{thm:relative-equilibria},
so $\zeta=k_s$. The lift along $\gamma$ is therefore
$\widehat Y(\gamma)=\tau X+k_sJX+(2k+c_0)\xi$, with $X=\gamma_s$.
Subtracting from \eqref{eq:relative-equilibrium-velocity} gives
$\gamma_t-\widehat Y(\gamma)=-v\gamma_s$, where $v=\tau-\tfrac12k^2$.
Arclength preservation gives $\tau_s=k\zeta=kk_s$ by \eqref{eq:sas-zc-a}, so $v_s=0$.
The resulting shift in arclength gives
\[
 \gamma(s,t)=\exp(t\widehat Y)\gamma_0(s-vt),\qquad k(s,t)=K(s-vt).
\]

\begin{proposition}[Traveling-wave reduction]
\label{prop:traveling-wave}
Let \(k(s,t)=K(s-vt)\).  Then \(k\) solves
\eqref{eq:sasakian-mkdv} on $G_\lambda$ if and only if there is a constant \(C\) such that
\begin{equation}
\label{eq:traveling-wave-ode}
 K''+\frac12K^3+(v+2\lambda)K=C.
\end{equation}
Equivalently, for another constant \(E_0\),
\begin{equation}
\label{eq:traveling-wave-first-integral}
 \frac12(K')^2+\frac18K^4+\frac{v+2\lambda}{2}K^2-CK=E_0.
\end{equation}
Thus periodic traveling waves are oscillatory trajectories in a quartic
potential.  The case \(v=0\) gives stationary curvature profiles.
\end{proposition}

\begin{proof}
Substitution in \eqref{eq:sasakian-mkdv} gives
\[
 -vK'=\left(K''+\frac12K^3+2\lambda K\right)'.
\]
One integration gives \eqref{eq:traveling-wave-ode}.  Multiplication by
\(K'\) and a second integration gives
\eqref{eq:traveling-wave-first-integral}.  The converse follows by
differentiating \eqref{eq:traveling-wave-ode}.
\end{proof}

We now normalize the two turning values; no nonconstant periodic profile
is lost.  Let \(K_-<K_+\) be the turning values of a nonconstant
oscillatory solution of \eqref{eq:traveling-wave-ode}.  Using the symmetry
\[
 (K,C)\longmapsto(-K,-C)
\]
if necessary, we may arrange that the turning value of larger absolute value
is positive.  Write it as \(b>0\), and write the other turning value as
\(br\), where \(-1\le r<1\).  With
\[
 u=bz,\qquad K(z)=b\,w(u),\qquad
 \nu=\frac{v+2\lambda}{b^2},\qquad
 c_\ast=\frac{C}{b^3},\qquad e_\ast=\frac{E_0}{b^4},
\]
the profile equation and first integral become
\begin{equation}
\label{eq:normalized-profile-ode}
 w_{uu}+\frac12w^3+\nu w=c_\ast,
 \qquad
 w_u^2=P_{r,\nu}(w),
\end{equation}
where the two turning values of \(w\) are exactly \(r\) and \(1\).  The
conditions \(P_{r,\nu}(r)=P_{r,\nu}(1)=0\) force
\begin{equation}
\label{eq:rotating-parameters}
 c_\ast(r,\nu)=\frac{(1+r)(1+r^2+4\nu)}8,
 \qquad
 e_\ast(r,\nu)=\frac18+\frac\nu2-c_\ast(r,\nu),
\end{equation}
and hence
\begin{equation}
\label{eq:rotating-polynomial}
 P_{r,\nu}(x)
 =2e_\ast-\frac14x^4-\nu x^2+2c_\ast x
 =\frac14(x-r)(1-x)Q_{r,\nu}(x),
\end{equation}
with
\begin{equation}
\label{eq:rotating-factorization}
 Q_{r,\nu}(x)
 =x^2+(1+r)x+1+r+r^2+4\nu.
\end{equation}
Conversely, every choice of \(b>0\), \(r\in[-1,1)\), and \(\nu\) for
which \(Q_{r,\nu}>0\) on \([r,1]\) gives a periodic traveling-wave profile on
\(G_\lambda\) by
\begin{equation}
\label{eq:denormalization}
 K(z)=b\,w(bz),\qquad
 v=b^2\nu-2\lambda,\qquad
 C=b^3c_\ast(r,\nu),\qquad E_0=b^4e_\ast(r,\nu).
\end{equation}
The factor \((x-r)(1-x)\) accounts for the two prescribed turning values.
A zero of \(Q_{r,\nu}\) at either endpoint would make it a multiple root of
\(P_{r,\nu}\), which a nonconstant profile cannot reach in finite time.
An interior zero would prevent the profile from crossing that value:
at a simple zero \(P_{r,\nu}\) changes sign, and at a multiple zero
the travel time diverges. Thus an oscillation between \(r\) and \(1\)
of finite period requires \(Q_{r,\nu}>0\) throughout \([r,1]\).
We call an admissible profile \emph{type~(a)} if the two roots of
\(Q_{r,\nu}\) form a nonreal conjugate pair, and \emph{type~(b)} if they
are real, counting a double root twice. The discriminant of \(Q_{r,\nu}\)
is \(-3r^2-2r-3-16\nu\). For \(r\le-1/3\) its vertex lies in \([r,1]\),
so positivity there forces type~(a). For \(r>-1/3\) the vertex lies to
the left of \(r\), and positivity on \([r,1]\) is equivalent to
\(Q_{r,\nu}(r)>0\). Thus type~(b) occurs exactly in the region
\begin{equation}
\label{eq:real-root-region}
 -\tfrac13<r<1,\qquad
 -\frac{3r^2+2r+1}{4}<\nu\le\nu_{\rm d}(r):=-\frac{3r^2+2r+3}{16}.
\end{equation}
Both real roots are then strictly below \(r\).
In Figure~\ref{fig:flat-closing-map}, type~(a) is the blue region
\(\nu>\nu_{\rm d}(r)\) for \(-1\le r<1\), and type~(b) is the tan
region to the right of the open diamond at \((-1/3,-1/6)\).

\begin{proposition}[General normalized profile]
\label{prop:general-cn-wave}
Let \(-1\le r<1\) and \(Q_{r,\nu}>0\) on \([r,1]\).
The two profile types have the following parametrizations.
\begin{enumerate}
\item[(a)] \emph{Conjugate roots.} Write
\begin{equation}
\label{eq:normalized-complex-roots}
 Q_{r,\nu}(x)=(x-a)^2+d^2,\qquad
 a=-\frac{1+r}{2},\qquad
 d^2=\frac{3r^2+2r+3}{4}+4\nu>0,
\end{equation}
\begin{equation}
\label{eq:general-cn-modulus-frequency}
 A_\pm=\sqrt{(x_\pm-a)^2+d^2},\quad x_+=1,\ x_-=r,\quad
 m=\frac{(1-r)^2-(A_+-A_-)^2}{4A_+A_-},\qquad
 \omega=\frac{\sqrt{A_+A_-}}2.
\end{equation}
Then \(0<m<1\), and after translating \(u\),
\begin{equation}
\label{eq:general-cn-profile}
 w(u)=
 \frac{A_-+rA_++(A_--rA_+)\cn(\omega u;m)}
 {A_++A_-+(A_--A_+)\cn(\omega u;m)},
 \qquad L_0=\frac{4\mathbb K(m)}{\omega}.
\end{equation}
\item[(b)] \emph{Real roots.}  If \(Q_{r,\nu}\) has real zeros
\(q_1\le q_2<r\), put
\begin{equation}
\label{eq:general-sn-modulus-frequency}
 m=\frac{(1-r)(q_2-q_1)}{(1-q_2)(r-q_1)},\qquad
 \omega=\frac{\sqrt{(1-q_2)(r-q_1)}}4.
\end{equation}
Then \(0\le m<1\), and after translating \(u\),
\begin{equation}
\label{eq:general-sn-profile}
 w(u)=\frac{(r-q_1)+q_1(1-r)\sn^2(\omega u;m)}
 {(r-q_1)+(1-r)\sn^2(\omega u;m)},
 \qquad L_0=\frac{2\mathbb K(m)}{\omega}.
\end{equation}
For a double root, \(m=0\) and the profile is elementary.
\end{enumerate}
In both cases \(K(s,t)=b\,w(b(s-vt))\), \(v=b^2\nu-2\lambda\), is a
traveling-wave solution on \(G_\lambda\) of spatial period \(L_0/b\), and
every nonconstant periodic profile arises in exactly one of the two forms.
\end{proposition}

\begin{proof}
In (a) the triangle inequalities $|A_+-A_-|<1-r<A_++A_-$ give $0<m<1$; in
(b), $(1-q_2)(r-q_1)-(1-r)(q_2-q_1)=(r-q_2)(1-q_1)>0$ gives $m<1$.  With
$y=\cn(\omega u;m)$, $(y')^2=\omega^2(1-y^2)(1-m+my^2)$, respectively
$y=\sn^2(\omega u;m)$, $(y')^2=4\omega^2y(1-y)(1-my)$, a direct
substitution gives $w_u^2=P_{r,\nu}(w)$, with $w=1$ at $y=1$, respectively
$y=0$, and $w=r$ at $y=-1$, respectively $y=1$.  The periods are those of
$\cn$ and $\sn^2$.
\end{proof}

Closure of the reconstructed curve will be governed by two integrals over
one curvature period.  Let \(w\) be the normalized periodic profile of period \(L_0\),
with the origin chosen so that \(w(0)=1\). Set
\begin{equation}
\label{eq:rotating-area-function}
  \Theta(r,\nu)=\int_0^{L_0}w\,du,
  \qquad
  \mathcal B(r,\nu)=\int_0^{L_0}\left(\frac12w^2+\nu\right)du.
\end{equation}
Thus \(\Theta\) is the total curvature of one period, that is, the turning of
the tangent, and \(\mathcal B\) will control the signed area swept in one
period (Remark~\ref{rem:noether-criterion} and
Lemma~\ref{lem:split-rotating-data}).  For \(r<1\), the substitution
\(x=w\) on the two half-periods gives
\begin{equation}
\label{eq:rotating-turning}
  \Theta(r,\nu)=2\int_r^1\frac{x\,dx}{\sqrt{P_{r,\nu}(x)}},
  \qquad
  \mathcal B(r,\nu)=2\int_r^1
  \frac{\frac12x^2+\nu}{\sqrt{P_{r,\nu}(x)}}\,dx,
\end{equation}
and the further substitution
\begin{equation}
\label{eq:x-r-phi}
  x_r(\phi)=\frac{1+r}{2}+\frac{1-r}{2}\cos\phi,
  \qquad 0\le\phi\le\pi,
\end{equation}
removes the endpoint singularities:
\begin{align}
  \Theta(r,\nu)
  &=4\int_0^\pi
    \frac{x_r(\phi)}{\sqrt{Q_{r,\nu}(x_r(\phi))}}\,d\phi,
  \label{eq:rotating-turning-regular}\\
  \mathcal B(r,\nu)
  &=4\int_0^\pi
    \frac{\frac12x_r(\phi)^2+\nu}
    {\sqrt{Q_{r,\nu}(x_r(\phi))}}\,d\phi.
  \label{eq:rotating-area-regular}
\end{align}
These regular forms also make sense at \(r=1\), and we use them to define
\(\Theta\) and \(\mathcal B\) there.

Since \(C=b^3c_\ast\) and \(8c_\ast=(1+r)(1+r^2+4\nu)\)
by \eqref{eq:rotating-parameters}, the admissible locus \(C=0\) in
Figure~\ref{fig:flat-closing-map} has two branches.
On the left edge \(r=-1\), positivity of \(Q_{-1,\nu}(x)=x^2+1+4\nu\)
requires \(\nu>-1/4\). Here \(w\) oscillates between \(-1\) and \(1\)
(the \(\cn\) profile of type~(a)).
On the parabola \(\nu=-(1+r^2)/4\), one has
\(Q_{r,\nu}(x)=(x+1)(x+r)\), so the admissible arc has \(0<r\le1\).
For \(0<r<1\) the roots \(-1,-r\) give the positive profile
\(w(u)=\dn(u/2;1-r^2)\) of type~(b).
Since \(P_{r,\nu}(x)=\tfrac14(1-x^2)(x^2-r^2)\), substituting
\(h=x^2\) in \eqref{eq:rotating-turning} gives
\[
 \Theta=4\int_r^1\frac{x\,dx}{\sqrt{(1-x^2)(x^2-r^2)}}
 =2\int_{r^2}^1\frac{dh}{\sqrt{(1-h)(h-r^2)}}=2\pi.
\]
The endpoint \((1,-1/2)\) is the small-amplitude limit, where
\eqref{eq:rotating-turning-regular} also gives \(\Theta=2\pi\).

The scalar profile depends on the Sasakian member only through
$v=b^2\nu-2\lambda$; the reconstruction in $\Sigma_{2\lambda}$ depends
essentially on $\lambda$.

\subsection{Space-form monodromy and Noether reconstruction}
\label{sec:space-form-reconstruction}

We now pass from the periodic scalar profile to the reconstructed curve on
its constant-curvature Reeb quotient.  Put
\[
  \kappa=2\lambda,
\]
and let \(\Sigma_\kappa\) denote the simply connected oriented surface of
Gaussian curvature \(\kappa\).

Let \(K\) be a smooth \(L\)-periodic function and \(\beta:\mathbb R\to
\Sigma_\kappa\) the unit-speed curve with signed geodesic curvature \(K\),
with Frenet frame \((T,N)\).  By the fundamental theorem of curves, the
unique \(\mathcal M_K\in\operatorname{Isom}^+(\Sigma_\kappa)\) carrying the
frame at \(s=0\) to that at \(s=L\) satisfies
\begin{equation}
\label{eq:space-form-monodromy}
 \beta(s+L)=\mathcal M_K\beta(s),\qquad
 T(s+L)=d\mathcal M_K\,T(s),
\end{equation}
so \(\beta\) closes smoothly after \(n\) curvature periods if and only if
\begin{equation}
\label{eq:space-form-base-closing}
 \mathcal M_K^n=\mathrm{id};
\end{equation}
Legendrian closure of the lift is then decided by
Theorem~\ref{thm:unified-holonomy}.  This is standard; for $S^3$ it is
\cite[Remark~4.4]{CIM2024}.

The traveling-wave equation determines the monodromy more rigidly because
its projected velocity is itself a Killing field.  This is the Killing
field along a constrained elastic curve, cf.\
\cite{LangerSinger1984,Heller2014}; for $\kappa>0$ it is the momentum of
\cite[(4.6)--(4.7)]{CIM2024}, whose eigenvalues correspond to
$\pm\sqrt{\mathcal I}$ in \eqref{eq:killing-casimir}.  We record it
uniformly in $\kappa$, since the split case below uses $(Y,\psi)$ and
$\mathcal I$ explicitly.

\begin{proposition}[Noether--Killing reconstruction on a space form]
\label{prop:space-form-noether-reconstruction}
Let \(K\) satisfy
\begin{equation}
\label{eq:space-form-profile-again}
  K_{ss}+\frac12K^3+(v+\kappa)K=C
\end{equation}
along a unit-speed curve \(\beta\subset\Sigma_\kappa\).  Define along
\(\beta\)
\begin{equation}
\label{eq:space-form-noether-field}
  Y=\left(\frac12K^2+v\right)T+K_sN,
  \qquad
  \psi=C-\kappa K.
\end{equation}
Then
\begin{equation}
\label{eq:killing-transport-along-curve}
  \nabla_TY=\psi N,
  \qquad
  \psi_s=-\kappa\langle Y,N\rangle.
\end{equation}
These are the Killing-transport equations on \(\Sigma_\kappa\).  Hence
\((Y,\psi)\) is the restriction of a unique global Killing field on the
simply connected space form.  The projected Hamiltonian velocity satisfies
\begin{equation}
\label{eq:space-form-relative-velocity}
  \frac12K^2T+K_sN=Y-vT,
\end{equation}
so, after removing the constant parameter translation, the projected
Hamiltonian motion is precisely the one-parameter isometry generated by
\(Y\).

If \(f_Y\) is the Killing potential normalized by
\(df_Y=-2\iota_YdA_\kappa\), then along the curve
\begin{equation}
\label{eq:killing-potential-curvature}
  f_Y\circ\beta=2K+c_0
\end{equation}
for a constant \(c_0\).  Thus the reconstruction is exactly the local
relative-equilibrium criterion of Theorem~\ref{thm:relative-equilibria}.
Finally,
\begin{equation}
\label{eq:killing-casimir}
  \mathcal I:=\psi^2+\kappa|Y|^2
  =C^2+\kappa(2E_0+v^2)
\end{equation}
is constant.  For a nonzero generator and \(\kappa<0\), the signs
\(\mathcal I>0\), \(\mathcal I=0\), and \(\mathcal I<0\) correspond,
respectively, to elliptic, parabolic, and hyperbolic one-parameter
subgroups; for \(\kappa=0\), the cases \(C=0\) and \(C\ne0\) become the
translational and rotational Euclidean branches.
\end{proposition}

\begin{proof}
Put
\[
  \tau_Y=\frac12K^2+v,\qquad \zeta_Y=K_s,
\]
so that \(Y=\tau_YT+\zeta_YN\).  Since \((\tau_Y)_s=K\zeta_Y\), the Frenet equations give
\[
  \nabla_TY=((\zeta_Y)_s+K\tau_Y)N
  =\left(K_{ss}+\frac12K^3+vK\right)N
  =(C-\kappa K)N=\psi N.
\]
Also \(\psi_s=-\kappa K_s=-\kappa\langle Y,N\rangle\), proving
\eqref{eq:killing-transport-along-curve}.  Killing transport on a simply
connected constant-curvature surface identifies the three-dimensional
space of initial data \((Y,\psi)\) at one point with the Lie algebra of
global Killing fields, so the pair extends uniquely.  Equation
\eqref{eq:space-form-relative-velocity} is immediate from the definition of
\(Y\).

With the orientation convention \(dA_\kappa(T,N)=1\),
\[
  \frac{d}{ds}(f_Y\circ\beta)
  =-2dA_\kappa(Y,T)=2K_s,
\]
which proves \eqref{eq:killing-potential-curvature}.  Finally,
\[
  \frac{d}{ds}\bigl(\psi^2+\kappa|Y|^2\bigr)
  =-2\kappa\psi K_s+2\kappa\psi\langle N,Y\rangle=0.
\]
The first integral \eqref{eq:traveling-wave-first-integral} gives
\[
  |Y|^2
  =K_s^2+\left(\frac12K^2+v\right)^2
  =2E_0+v^2-\kappa K^2+2CK,
\]
and substitution of \(\psi=C-\kappa K\) yields
\eqref{eq:killing-casimir}.  The final classification is the standard
classification of one-parameter subgroups of the isometry group of a
space form; in the flat case it is also verified directly below.
\end{proof}

\subsection{Closed relative equilibria in the Heisenberg group}
\label{sec:rotating-torus}

We next specialize the closing data to the flat member
\(\kappa=0\).  In this case one may choose Euclidean coordinates and write
\[
  \theta(s)=\theta_0+\int_0^sK(s')\,ds',
  \qquad
  \beta(s)=\beta_0+\int_0^s(\cos\theta(s'),\sin\theta(s'))\,ds'.
\]
For one curvature period \(L\), smooth base closure is therefore equivalent
to
\begin{equation}
\label{eq:planar-closing-conditions}
  \int_0^L K\,ds=2\pi n,
  \qquad
  \int_0^L e^{i\theta(s)}\,ds=0,
  \qquad n\in\mathbb Z,
\end{equation}
and the horizontal lift to \(H^1\) closes precisely when, in addition,
\begin{equation}
\label{eq:curvature-zero-area}
  \oint_\beta(x\,dy-y\,dx)=0.
\end{equation}
These are the \(\kappa=0\) form of \eqref{eq:space-form-base-closing}
and Theorem~\ref{thm:unified-holonomy}. The second integral in
\eqref{eq:planar-closing-conditions} is the displacement
$\beta(L)-\beta(0)$ in $\mathbb C$. To express the closing data in normalized
variables, set $u=bs$. Then $K=bw(u)$, $v=b^2\nu$, $C=b^3c_\ast(r,\nu)$,
and $L=L_0/b$. Since $(K^2/2+v)\,ds=b(w^2/2+\nu)\,du$, the definitions
in \eqref{eq:rotating-area-function} give
\begin{equation}
\label{eq:scaled-closing-data}
  \int_0^L K\,ds=\Theta,
  \qquad
  \int_0^L\left(\frac12K^2+v\right)ds=b\mathcal B.
\end{equation}
\begin{remark}[Flat closing data]
\label{rem:noether-criterion}
The first identity in \eqref{eq:scaled-closing-data} gives the change in
tangent angle. The second determines the displacement for $C=0$ by
\eqref{eq:translation-monodromy-B} and the signed area contribution for
$C\ne0$ by \eqref{eq:rotational-area-B}. To see this, put
\begin{equation}
\label{eq:euclidean-noether-vector}
  \mathcal J
  =K_sT-\left(\frac12K^2+v\right)N.
\end{equation}
Then
\begin{equation}
\label{eq:noether-velocity}
  \frac12K^2T+K_sN=J\mathcal J-vT.
\end{equation}

In the branch \(C=0\), Proposition~\ref{prop:space-form-noether-reconstruction}
identifies the projected motion as a translation. Indeed, differentiating
\eqref{eq:euclidean-noether-vector} with $T_s=KN$ and $N_s=-KT$ cancels
the normal terms and gives $\mathcal J_s=(K_{ss}+K^3/2+vK)T=CT$
by \eqref{eq:traveling-wave-ode} at $\lambda=0$. Thus $\mathcal J$ is
constant when $C=0$. In the fixed orthogonal basis
$(\mathcal J,J\mathcal J)$, the tangent has scalar products
$\langle T,\mathcal J\rangle=K_s$ and
$\langle T,J\mathcal J\rangle=K^2/2+v$. Integrating $T=\beta_s$
over one period and using $K(L)-K(0)=0$ gives
\begin{equation}
\label{eq:translation-monodromy-B}
  \beta(L)-\beta(0)
  =\frac{1}{|\mathcal J|^2}
   \left(\int_0^L\left(\frac12K^2+v\right)ds\right)J\mathcal J.
\end{equation}
For a nonconstant profile $|\mathcal J|>0$.  Hence, once the tangent
monodromy is trivial, position closure is
equivalent to \(\mathcal B(r,\nu)=0\).  The Legendrian zero-area condition
\eqref{eq:curvature-zero-area} remains an additional requirement; for the
symmetric figure-eight of Theorem~\ref{thm:closed-elliptic} it follows
from half-period reflection.

The positive $\dn$ branch of the case $C=0$ has trivial tangent
monodromy but, by Proposition~\ref{prop:real-root-no-closure}, nonzero
translational displacement.  Together with Theorem~\ref{thm:closed-elliptic},
this accounts for the nonconstant closed flat $C=0$ profiles
(Corollary~\ref{cor:translational-figure-eight}), in agreement
with the closed-elastica classification \cite{LangerSinger1984}.

For $C\ne0$, the same identity gives $(\beta-\mathcal J/C)_s=0$ since
$\beta_s=T$. Write $o=\beta-\mathcal J/C$. The Killing field along $\beta$
is $J\mathcal J=CJ(\beta-o)$ by \eqref{eq:noether-velocity}, so $o$ is the
rotation center. Taking $o$ as the origin gives
\begin{equation}
\label{eq:general-rotating-reconstruction}
  \beta=\frac{\mathcal J}{C}
  =\frac1C\left\{
  K_sT-\left(\frac12K^2+v\right)N\right\}.
\end{equation}
Over one curvature period \(\beta\) is rotated through the angle
\(\Theta(r,\nu)\); therefore
\begin{equation}
\label{eq:rotational-angle-closing}
  \Theta(r,\nu)=\frac{2\pi p}{n},\qquad (p,n)=1,
\end{equation}
closes both position and tangent after \(n\) curvature periods.  Moreover,
\begin{equation}
\label{eq:general-rotating-area}
  \det(\beta,\beta_s)
  =\frac1C\left(\frac12K^2+v\right),
\end{equation}
so twice the signed-area contribution of one curvature period is
\begin{equation}
\label{eq:rotational-area-B}
  \int_0^L\det(\beta,\beta_s)\,ds
  =\frac{\mathcal B(r,\nu)}{b^2c_\ast(r,\nu)}.
\end{equation}
Thus the Legendrian closing condition is exactly \(\mathcal B(r,\nu)=0\).  Finally,
\begin{equation}
\label{eq:general-rotating-radius}
  C^2|\beta|^2
  =K_s^2+\left(\frac12K^2+v\right)^2
  =2E_0+v^2+2CK.
\end{equation}
Hence in the rotational branch \(\Theta\) and \(\mathcal B\) are precisely the
base-monodromy and Reeb-holonomy closing data used below.
\end{remark}

\subsubsection*{The branch \(C=0\): the closed figure-eight}

\begin{theorem}[Closed elliptic relative equilibrium]
\label{thm:closed-elliptic}
There is a unique \(m_*\in(0,1)\) satisfying
\begin{equation}
\label{eq:elliptic-closing-modulus}
 2\mathbb E(m_*)=\mathbb K(m_*),
 \qquad m_*=0.8261147659\ldots.
\end{equation}
For every \(b>0\), set
\begin{equation}
\label{eq:cn-traveling-wave}
 k(s,t)=2b\sqrt{m_*}\,
 \cn\!\left(b(s-v_0t);m_*\right),
 \qquad
 v_0=(1-2m_*)b^2.
\end{equation}
Then \(k\) is an \(\ell\)-periodic solution of the Sasakian mKdV equation
on \(G_0=H^1\), where \(\ell=4\mathbb K(m_*)/b\).  Its planar
reconstruction is a closed unit-speed figure-eight elastica with zero
signed area.  Consequently it has a closed Legendrian lift, unique up to
Reeb translation, and the Hamiltonian flow reconstructs it as a relative
equilibrium modulo pseudo-Hermitian rigid motions and parameter
translations.  
\end{theorem}

\begin{figure}[!htbp]
 \centering
 \includegraphics[width=\linewidth]{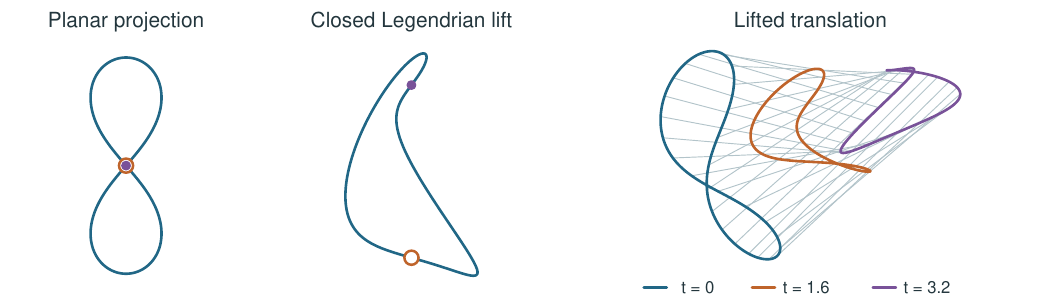}
 \caption{The translating figure-eight of Theorem~\ref{thm:closed-elliptic}
 with $b=1$. The points at $u=\mathbb K(m_*)$ (open circle) and
 $u=3\mathbb K(m_*)$ (solid dot) coincide in the planar projection
 (left) and have distinct heights on the Legendrian lift (middle).
 At right, the lift is shown at three times under the left translation
 $\Phi_t(x,y,z)=(x+t,y,z-ty)$. The thin lines trace its point orbits.}
 \label{fig:flat-eight-motion}
\end{figure}

\begin{proof}
The Jacobi identity
\[
 \cn''(u;m)=(2m-1)\cn(u;m)-2m\cn^3(u;m)
\]
shows that
\[
 K(z)=2b\sqrt m\,\cn(bz;m)
\]
satisfies
\[
 K''+\frac12K^3+(1-2m)b^2K=0.
\]
This is the symmetric specialization \(C=0\) of
Proposition~\ref{prop:general-cn-wave}.  Since
\(v_0=(1-2m)b^2\) and \(\lambda=0\),
Proposition~\ref{prop:traveling-wave}
proves the mKdV assertion, and the period is \(4\mathbb K(m)/b\).

The ratio \(\mathbb E(m)/\mathbb K(m)\) decreases strictly from \(1\) to
\(0\) as \(m\) increases from \(0\) to \(1\).  Indeed, the standard
derivative identities
\[
 \mathbb E'=\frac{\mathbb E-\mathbb K}{2m},
 \qquad
 \mathbb K'=\frac{\mathbb E}{2m(1-m)}-\frac{\mathbb K}{2m}
\]
give, with \(\varrho=\mathbb E/\mathbb K\),
\[
 \varrho'
 =\frac{2\varrho-1-\varrho^2/(1-m)}{2m}<0,
\]
because
\[
 \varrho^2-2(1-m)\varrho+(1-m)>0
\]
has negative discriminant for \(0<m<1\).  This proves existence and
uniqueness of \(m_*\).

Fix a time and write \(u=b(s-v_0t)\).  After a planar rigid motion, choose
\[
 \theta(u)=2\arcsin\!\bigl(\sqrt m\,\sn(u;m)\bigr).
\]
Since \(\dn(u;m)>0\) for \(0<m<1\), differentiation gives
\[
 \theta_s=2b\sqrt m\,\cn(u;m)=k.
\]
Consequently
\[
 (\cos\theta,\sin\theta)
 =
 \bigl(1-2m\sn^2u,\;2\sqrt m\,\sn u\,\dn u\bigr).
\]
With the Jacobi epsilon function
$\mathcal E(u;m)=\int_0^u\dn^2(w;m)\,dw$ and $1-2m\sn^2u=2\dn^2u-1$, a
unit-speed planar reconstruction is therefore
\begin{equation}
\label{eq:explicit-elastica}
 x(u)=\frac{2\mathcal E(u;m)-u}{b},
 \qquad
 y(u)=-\frac{2\sqrt m}{b}\cn(u;m).
\end{equation}
Using
\[
 \mathcal E(u+4\mathbb K;m)
 =\mathcal E(u;m)+4\mathbb E(m)
\]
and the \(4\mathbb K\)-periodicity of \(\cn\), one finds
\[
 x(u+4\mathbb K)-x(u)
 =\frac4b\bigl(2\mathbb E-\mathbb K\bigr),
 \qquad
 y(u+4\mathbb K)=y(u).
\]
Thus \eqref{eq:elliptic-closing-modulus} is exactly the position-closing
condition.  The tangent is already \(4\mathbb K\)-periodic and its total
turning is zero, so both conditions in
\eqref{eq:planar-closing-conditions} hold.

It remains to verify Legendrian closure.  The half-period identities give
\[
 x(u+2\mathbb K)=x(u),\qquad
 y(u+2\mathbb K)=-y(u)
\]
when \(2\mathbb E=\mathbb K\).  Hence the second half of the curve is the
reflection of the first half across the \(x\)-axis, with the same
parameter orientation.  Reflection reverses signed area, so the two
contributions to \(\oint(x\,dy-y\,dx)\) cancel, and
Theorem~\ref{thm:unified-holonomy} gives a closed Legendrian lift.
Moreover, \(u=\mathbb K\) and \(u=3\mathbb K\) both give the point
\((x,y)=(0,0)\), and the two intervening arcs are reflected lobes.  Thus
the closed planar curve is the classical figure-eight elastica
\cite{LangerSinger1984}.

Finally, this is the flat \(C=0\) specialization of
Proposition~\ref{prop:space-form-noether-reconstruction}, as made explicit
in Remark~\ref{rem:noether-criterion}.  After the parameter translation is
removed, the projected curve therefore evolves by translation, and the lifted translation generates the Hamiltonian trajectory
of the closed Legendrian curve on $\mathcal Q$. Figure~\ref{fig:flat-eight-motion} shows
the projection, its lift, and the lifted translation.
\end{proof}

Constant curvature gives no closed Legendrian curve in $H^1$: for $k\ne0$
the planar curve is a circle with nonzero signed area, and for $k=0$ it is a
line.

\begin{proposition}[No closed real-root profiles in $H^1$]
\label{prop:real-root-no-closure}
In the region \eqref{eq:real-root-region},
\[
 \mathcal B(r,\nu)\le\mathcal B(r,\nu_{\rm d}(r))
 =-\frac{\pi(1-r)^2}{2\sqrt{(3r+1)(r+3)}}<0 .
\]
Consequently, the curvature profile of every closed Legendrian relative
equilibrium in the Heisenberg group $H^1$ is of type~\textup{(a)} in
Proposition~\ref{prop:general-cn-wave}, so the two roots of $Q_{r,\nu}$
form a nonreal conjugate pair.
\end{proposition}

\begin{proof}
For $r>-\frac13$ and $x\ge r$, the $\nu$-derivative of the integrand of
\eqref{eq:rotating-area-regular} has numerator
$(1+r)x+1+r+r^2+2\nu\ge\frac12(1+r)^2>0$, so $\mathcal B$ increases in
$\nu$.  At $\nu=\nu_{\rm d}$, $\sqrt{Q_{r,\nu}}=y:=x+\frac{1+r}2$ and
$\frac12x^2+\nu_{\rm d}=\frac12y^2-\frac{1+r}2y-\frac{(1-r)^2}{16}$; with
$y=(1+r)+\frac{1-r}2\cos\phi$ the integral is elementary and gives the
value shown.  Closure requires $\mathcal B=0$: by
Remark~\ref{rem:noether-criterion} for $C=0$, and by
\eqref{eq:rotational-area-B} for $C\ne0$.
\end{proof}

\subsubsection*{The branch \(C\ne0\): rotating relative equilibria}

We now turn to the nonzero integration-constant branch.  The normalized
family \(P_{r,\nu}\), \(Q_{r,\nu}\), \(c_\ast(r,\nu)\), and \(e_\ast(r,\nu)\) was
introduced before Proposition~\ref{prop:general-cn-wave}, and the closing
quantities \(\Theta\) and \(\mathcal B\) were defined in
\eqref{eq:rotating-area-function}.
The remaining task is therefore purely analytic: find admissible pairs
\((r,\nu)\) for which
\[
 \Theta(r,\nu)=\Theta_0,\qquad \mathcal B(r,\nu)=0,
\]
with \(0<\Theta_0<2\pi\), and show that the resulting pair lies in the
rotational branch \(c_\ast(r,\nu)\ne0\).  We first solve the area equation for $\nu$ at each $r$, then follow the
turning angle along the resulting graph.

\begin{equation}
\label{eq:nu-min}
 \nu_{\min}(r)=
 \begin{cases}
 -\dfrac{3r^2+2r+3}{16}, & -1\le r\le-\tfrac13,\\[2mm]
 -\dfrac{3r^2+2r+1}{4}, & -\tfrac13\le r\le1,
 \end{cases}
\end{equation}
a continuous function (both expressions equal \(-\tfrac16\) at
\(r=-\tfrac13\)) with \(\nu_{\min}\le-\tfrac16<0\).

\begin{lemma}[Admissible region]
\label{lem:admissible-region}
For \(r\in[-1,1]\), one has \(Q_{r,\nu}>0\) on \([r,1]\) if and only if
\(\nu>\nu_{\min}(r)\).  If \(\nu=\nu_{\min}(r)\), then \(Q_{r,\nu}\)
vanishes on \([r,1]\) exactly at
\[
 x_d(r)=
 \begin{cases}
 -\dfrac{1+r}{2}, & -1\le r\le-\tfrac13,\\[2mm]
 r, & -\tfrac13\le r\le1,
 \end{cases}
\]
and
\[
 \frac12x_d(r)^2+\nu_{\min}(r)=
 \begin{cases}
 -\dfrac{(1-r)^2}{16}, & -1\le r\le-\tfrac13,\\[2mm]
 -\dfrac{(1+r)^2}{4}, & -\tfrac13\le r\le1.
 \end{cases}
\]
In particular \(\tfrac12x_d(r)^2+\nu_{\min}(r)<0\) for every
\(r\in[-1,1]\).
\end{lemma}

\begin{proof}
Write \(Q_{r,\nu}(x)=\bigl(x+\tfrac{1+r}2\bigr)^2
+\tfrac{3r^2+2r+3}{4}+4\nu\).  The vertex \(-\tfrac{1+r}2\) lies in
\([r,1]\) if and only if \(r\le-\tfrac13\); then the minimum of
\(Q_{r,\nu}\) on \([r,1]\) is \(\tfrac{3r^2+2r+3}{4}+4\nu\), attained only
at the vertex.  If \(r\geq-\tfrac13\), \(Q_{r,\nu}\) is strictly
increasing on \([r,1]\), and its minimum there is
\(Q_{r,\nu}(r)=3r^2+2r+1+4\nu\), attained only at \(x=r\).  The last
formulas follow by substitution.
\end{proof}

\begin{figure}[!t]
 \centering
 \begin{tikzpicture}
  \node[anchor=south west,inner sep=0] (img)
   {\includegraphics[width=0.94\linewidth]{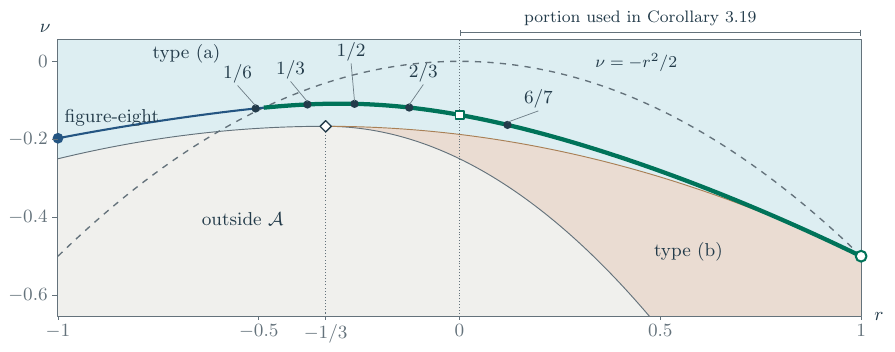}};
  \begin{scope}[x={(img.south east)},y={(img.north west)},
                font=\fontsize{9.4pt}{11pt}\selectfont]
   \fill[white] (0.583,0.918) rectangle (0.852,0.976);
   \node[anchor=base west,inner sep=0] at (0.5857,0.9330)
    {portion used in Corollary~\ref{cor:embedded-rotation-numbers}};
  \end{scope}
 \end{tikzpicture}
 \caption{Closing parameters in $H^1$. The solid graph is $\mathcal B=0$;
 its left endpoint is the figure-eight and the open right endpoint is the
 small-amplitude limit. The diamond at $(-1/3,-1/6)$ marks the left tip
 of the type~(b) region in \eqref{eq:real-root-region}. Dots labeled $p/n$ mark sample rotating curves
 that close after $n$ curvature periods. Below the dashed parabola,
 closing pairs on the thick part of the graph sweep embedded tori by
 Theorem~\ref{thm:rotating-torus}. Corollary~\ref{cor:embedded-rotation-numbers}
 uses the portion $0\le r<1$, starting at the square where
 $\Theta=\Theta^*$. The graph and marked pairs are plotted numerically.}
 \label{fig:flat-closing-map}
\end{figure}

Let
\[
 \mathcal{A}=\{(r,\nu):-1\le r\le1,\ \nu>\nu_{\min}(r)\}.
\]
Figure~\ref{fig:flat-closing-map} shows $\mathcal{A}$ as the union of the blue and tan regions
strictly above the lower boundary $\nu=\nu_{\min}(r)$, with vertical edges at $r=\pm1$.
For $r<1$, this region gives oscillations with simple turning values
$r$ and $1$; the edge $r=1$ is their small-amplitude limit.
Positivity only on the open interval would also allow multiple endpoint
roots, for which the period need not be finite.  The region $\mathcal{A}$ contains both cases of
Proposition~\ref{prop:general-cn-wave}; by
Proposition~\ref{prop:real-root-no-closure}, only type~(a) can yield a closed Legendrian lift.  The functions
\(\Theta\) and \(\mathcal B\) defined in
\eqref{eq:rotating-turning-regular} and
\eqref{eq:rotating-area-regular} are continuous on \(\mathcal{A}\), including
the edges \(r=\pm1\).

\begin{lemma}[Boundary behavior]
\label{lem:rotating-boundary}
\mbox{}
\begin{enumerate}
\item\label{it:left-edge} \(\Theta(-1,\nu)=0\) for every
\(\nu>-\tfrac14=\nu_{\min}(-1)\).
\item\label{it:right-edge} For \(\nu>-\tfrac32=\nu_{\min}(1)\),
\[
 \Theta(1,\nu)=\frac{4\pi}{\sqrt{6+4\nu}},
 \qquad
 \mathcal B(1,\nu)=\frac{4\pi\bigl(\frac12+\nu\bigr)}{\sqrt{6+4\nu}}.
\]
\item\label{it:top-edge} If \(\tfrac12x^2+\nu\geq0\) on \([r,1]\) and
\(\tfrac12+\nu>0\), then \(\mathcal B(r,\nu)>0\).
\item\label{it:bottom-edge} For each fixed \(r\in[-1,1]\),
\(\mathcal B(r,\nu)\to-\infty\) as \(\nu\downarrow\nu_{\min}(r)\).
\end{enumerate}
\end{lemma}

\begin{proof}
(\ref{it:left-edge}) For \(r=-1\) one has \(x_{-1}(\phi)=\cos\phi\) and
\(Q_{-1,\nu}(x)=x^2+1+4\nu\), so the integrand of
\eqref{eq:rotating-turning-regular} is odd under \(\phi\mapsto\pi-\phi\).

(\ref{it:right-edge}) For \(r=1\), \(x_1\equiv1\) and \(Q_{1,\nu}(1)=6+4\nu\).

(\ref{it:top-edge}) The integrand of \eqref{eq:rotating-area-regular} is
then nonnegative, and it is positive near \(\phi=0\), where
\(x_r(\phi)\) is close to \(1\).

(\ref{it:bottom-edge}) For \(r=1\) the limit follows from (\ref{it:right-edge}).
Fix \(r<1\) and put \(\nu_0=\nu_{\min}(r)\).
By Lemma~\ref{lem:admissible-region}, the function
\(q(\phi)=Q_{r,\nu_0}(x_r(\phi))\) has a unique zero \(\phi_d\) on
\([0,\pi]\), where \(x_r(\phi_d)=x_d(r)\). Since \(x_r(\phi)\in[r,1]\)
for every real \(\phi\), the smooth function \(q\) is nonnegative on
\(\mathbb R\). Thus \(q'(\phi_d)=0\) even if \(\phi_d=\pi\), and
\(q(\phi)=O((\phi-\phi_d)^2)\).
On a small interval about \(\phi_d\) in \([0,\pi]\), the numerator in
\eqref{eq:rotating-area-regular} stays below a negative constant for
\(\nu\) close to \(\nu_0\), since \(x_d(r)^2/2+\nu_0<0\).
As \(\nu\downarrow\nu_0\), the reciprocal denominator
\((q(\phi)+4(\nu-\nu_0))^{-1/2}\) increases to \(q(\phi)^{-1/2}\),
whose integral on that interval diverges. By monotone convergence the
interval contributes a quantity tending to \(-\infty\), while the
integral over its complement stays bounded.
\end{proof}

\begin{lemma}[The zero-area graph]
\label{lem:area-graph}
For each $r\in[-1,1]$ there is exactly one $\nu=\nu_B(r)$ with
$(r,\nu)\in\mathcal{A}$ and $\mathcal B(r,\nu)=0$. Its graph is real analytic,
with analytic extension at both endpoints. It joins the figure-eight
point $\bigl(-1,(1-2m_*)/(4m_*)\bigr)$ to $(1,-1/2)$.
On $-1<r<1$ the rotation coefficient $c_\ast(r,\nu_B(r))$ is positive.
\end{lemma}

\begin{proof}
With $x=x_r(\phi)$ fixed in \eqref{eq:rotating-area-regular}, differentiation
of the integrand gives
\[
 \partial_\nu\frac{x^2/2+\nu}{\sqrt{Q_{r,\nu}(x)}}
 =\frac{(1+r)x+1+r+r^2+2\nu}{Q_{r,\nu}(x)^{3/2}} .
\]
Since $x\ge r$, $1+r\ge0$ and $\nu>\nu_{\min}(r)$, the numerator exceeds
\[
 1+2r+2r^2+2\nu_{\min}(r)=
 \begin{cases}
 (13r^2+14r+5)/8,& -1\le r\le-1/3,\\
 (1+r)^2/2,&-1/3\le r\le1.
 \end{cases}
\]
Both bounds are positive: the quadratic in the first has discriminant
$-64$. Thus $\partial_\nu\mathcal B>0$ throughout $\mathcal{A}$.
For fixed $r$, parts (\ref{it:bottom-edge}) and (\ref{it:top-edge}) of
Lemma~\ref{lem:rotating-boundary} give $\mathcal B(r,\nu)\to-\infty$
as $\nu\downarrow\nu_{\min}(r)$ and $\mathcal B(r,0)>0$, respectively.
This proves existence and uniqueness. Near each zero the radicand in
\eqref{eq:rotating-area-regular} stays uniformly positive for
$0\le\phi\le\pi$. The regular integral is real analytic in $(r,\nu)$
there, even across $r=\pm1$ using the same formula. The analytic implicit
function theorem gives the asserted graph. Its value at $r=1$ follows
from Lemma~\ref{lem:rotating-boundary}(\ref{it:right-edge}); at $r=-1$,
normalizing the upper turning value of the figure-eight in
Theorem~\ref{thm:closed-elliptic} to one gives the stated value.

If $-1<r<1$ and $c_\ast=0$, the factorization
$Q_{r,\nu}(x)=(x+1)(x+r)$ forces $r>0$ by admissibility.
Proposition~\ref{prop:real-root-no-closure} then gives $\mathcal B<0$.
Thus $c_\ast$ cannot vanish on the interior graph. At $r=0$ it equals
$(1+4\nu_B(0))/8>0$ by admissibility, which fixes its sign throughout.
\end{proof}

\begin{corollary}[Closed translational relative equilibria]
\label{cor:translational-figure-eight}
Up to scaling, rigid motions of $H^1$, and Reeb translations, the closed
Legendrian lift of the figure-eight of Theorem~\ref{thm:closed-elliptic}
is the only closed relative equilibrium in $H^1$ with nonconstant
curvature whose generator is the lift of a translation.
\end{corollary}

\begin{proof}
By Remark~\ref{rem:noether-criterion}, the generator is the lift of a
translation exactly when $C=0$, that is, $c_\ast(r,\nu)=0$.  If $r>-1$,
this forces $1+r^2+4\nu=0$; admissibility then gives $r>0$, and the
profile is the positive $\dn$ profile
of type~(b), which does not close by
Proposition~\ref{prop:real-root-no-closure}.  Hence $r=-1$, and
$\Theta(-1,\nu)=0$ by Lemma~\ref{lem:rotating-boundary}.  By
\eqref{eq:translation-monodromy-B} each curvature period then adds the same
displacement, which vanishes only if $\mathcal B(-1,\nu)=0$.  By
Lemma~\ref{lem:area-graph} this fixes $\nu=(1-2m_*)/(4m_*)$, the
figure-eight of Theorem~\ref{thm:closed-elliptic}, whose lift is unique up
to Reeb translation.
\end{proof}

\begin{theorem}[Prescribed turning angles on the zero-area graph]
\label{thm:closing-pairs}
For every \(\Theta_0\in(0,2\pi)\) there is \((r,\nu)\in\mathcal{A}\) with
\(-1<r<1\) such that
\begin{equation}
\label{eq:rotating-closing-equations}
 \mathcal B(r,\nu)=0,\qquad \Theta(r,\nu)=\Theta_0,\qquad c_\ast(r,\nu)\neq0.
\end{equation}
Moreover every such pair satisfies \(\tfrac12+\nu>0\) and \(\nu<0\).
\end{theorem}

\begin{proof}
For each $r\in[-1,1]$, let $\nu_B(r)$ be the unique admissible value
for which $\mathcal B(r,\nu_B(r))=0$, as in Lemma~\ref{lem:area-graph}.
Along this zero-area graph, the turning angle is
$\vartheta(r)=\Theta(r,\nu_B(r))$. Since $\nu_B$ is continuous,
\eqref{eq:rotating-turning-regular} gives continuity of $\vartheta$ on $[-1,1]$. Lemma~\ref{lem:rotating-boundary} gives
$\vartheta(-1)=0$ and $\vartheta(1)=\Theta(1,-1/2)=2\pi$.
The intermediate value theorem yields $r\in(-1,1)$ with
$\vartheta(r)=\Theta_0$, and Lemma~\ref{lem:area-graph} gives $c_\ast>0$.
Finally the integrand in \eqref{eq:rotating-area-regular} must change
sign for every interior zero of $\mathcal B$. Since $|x|\le1$, this
implies $\nu<0<1/2+\nu$.
\end{proof}

On the rotating branch in $H^1$, we call an admissible pair $(r,\nu)$
a \emph{closing pair} if $\mathcal B(r,\nu)=0$ and
$\Theta(r,\nu)/(2\pi)\in\mathbb Q$.

\begin{remark}[Elliptic-integral form of $\Theta$ and $\mathcal B$]
\label{rem:closing-pair-elliptic}
By Proposition~\ref{prop:real-root-no-closure} every closing pair is of
type (a).  If $-1<r<1$ and $Q_{r,\nu}$ has no real zeros, Proposition~\ref{prop:general-cn-wave}
gives $w=\mathsf P+\mathsf Q/(1+\chi\cn(\omega u;m))$ with
$\chi=(A_--A_+)/(A_++A_-)$, $\mathsf P=(A_--rA_+)/(A_--A_+)$ and
$\mathsf Q=-2(1-r)A_+A_-/((A_--A_+)(A_++A_-))$.  Hence
\[
 \begin{gathered}
 \Theta=\frac{4\mathbb K\mathsf P+\mathsf QI_1}{\omega},\qquad
 \mathcal B=\frac{4\mathbb K(\frac12\mathsf P^2+\nu)+\mathsf P\mathsf QI_1
 +\frac12\mathsf Q^2I_2}{\omega},\\
 I_1=\frac{4\boldsymbol\Pi(-\frac{\chi^2}{1-\chi^2};m)}{1-\chi^2},\quad
 I_2=\frac{d(\chi I_1)}{d\chi}.
 \end{gathered}
\]
Both closing functions are combinations of complete elliptic integrals
of the three kinds.  These formulas are used only for numerical evaluation.
\end{remark}

\begin{theorem}[Rotating closed Legendrian curves]
\label{thm:rotating-family}
Let \(p/n\in(0,1)\) be a reduced fraction. Taking \(\Theta_0=2\pi p/n\)
in Theorem~\ref{thm:closing-pairs} gives a closing pair \((r,\nu)\).  Put \(c_\ast=c_\ast(r,\nu)\ne0\) and \(e_\ast=e_\ast(r,\nu)\).  Let $K$ be the normalized periodic oscillator with $K(0)=1$,
$K_s(0)=0$, and $K_s^2=P_{r,\nu}(K)$, of least period
$L=2\int_r^1P_{r,\nu}(x)^{-1/2}dx$.  Its traveling-wave curvature
$k(s,t)=K(s-\nu t)$ reconstructs a closed Legendrian curve
\(\gamma_0\subset H^1\) of period \(nL\), and the Hamiltonian vector field
of horizontal length moves it as a rotating relative equilibrium with
angular velocity \(c_\ast\).  There are infinitely many geometrically distinct closed rotating relative
equilibria.  The scale $b>0$ is restored by $K_b(s)=bK(bs)$, with speed
$b^2\nu$ and angular velocity $b^3c_\ast$.
\end{theorem}

\begin{proof}
Apply Remark~\ref{rem:noether-criterion} with \(b=1\), \(v=\nu\), and
\(C=c_\ast\).  Thus \(\beta\) is given by
\eqref{eq:general-rotating-reconstruction}, \(\beta_s=T\), and the
projected evolution is, modulo parameter translation, rotation with
angular velocity \(c_\ast\).
Since \(K\) is \(L\)-periodic with \(\int_0^LK\,ds=\Theta_0\), the frame
satisfies \(T(s+L)=\operatorname{Rot}_{\Theta_0}T(s)\), \(N(s+L)=\operatorname{Rot}_{\Theta_0}N(s)\), where
\(\operatorname{Rot}_\alpha\) denotes rotation through \(\alpha\) in the \(xy\)-plane.
Consequently \(\beta(s+L)=\operatorname{Rot}_{\Theta_0}\beta(s)\), and \(\beta\) closes after
\(n\) curvature periods.

By \eqref{eq:rotational-area-B}, twice the signed area swept over one
curvature period is \(\mathcal B(r,\nu)/c_\ast=0\).  In the coordinates
\eqref{eq:heisenberg-model}, the horizontal lift satisfies
\begin{equation}
\label{eq:rotating-z-lift}
 z_s=-\det(\beta,\beta_s)=-\frac1{c_\ast}\left(\frac12K^2+\nu\right),
\end{equation}
which is \(L\)-periodic with zero integral over one period.  Hence
\(\gamma_0=(\beta,z)\) is a closed Legendrian curve of period \(nL\), with
\(z(s+L)=z(s)\).

Combining rotation about the $z$-axis with the arclength shift $s\mapsto s-\nu t$, we obtain
\begin{equation}
\label{eq:rotating-flow-map}
 F(s,t)=\bigl(\operatorname{Rot}_{c_\ast t}\beta(s-\nu t),z(s-\nu t)\bigr).
\end{equation}
Equation~\eqref{eq:noether-velocity} shows at once that the planar
component of \(F_t\) is the horizontal part of
\eqref{eq:sasakian-curve-flow}.  Since planar rotations preserve the
contact form $\eta=dz+x\,dy-y\,dx$ of \eqref{eq:heisenberg-model},
\eqref{eq:general-rotating-radius} gives
\[
 \eta(F_t)=c_\ast|\beta|^2
 =2K+c_0,
 \qquad
 c_0=\frac{2e_\ast+\nu^2}{c_\ast}.
\]
Therefore
\[
 F_t=\frac12K^2X+K_sJX+(2K+c_0)\xi.
\]
The last summand \(c_0\xi\) is a constant Reeb translation and represents
zero on \(\mathcal Q\).  Thus $F$ represents the Hamiltonian trajectory of
Theorem~\ref{thm:sasakian-family-mkdv}. For this value of $c_0$,
\eqref{eq:relative-equilibrium-velocity} gives rotation about the $z$-axis up
to the arclength shift, with Killing potential $f_Y=c_\ast|\beta|^2$. The
rotation period is \(2\pi/|c_\ast|\).  Different fractions \(p/n\) give
different turning numbers \(p\) of \(\beta\) over its period, or different
numbers \(n\) of curvature periods, and hence geometrically distinct
curves.
\end{proof}

For the embeddedness of the swept torus we need the sign condition
\begin{equation}
\label{eq:torus-condition}
 \frac12r^2+\nu<0.
\end{equation}

\begin{theorem}[Embedded swept torus]
\label{thm:rotating-torus}
In Theorem~\ref{thm:rotating-family}, assume in addition
\eqref{eq:torus-condition}.  Then, for the representative of the
Hamiltonian vector field with constant Reeb term \(c_0\xi\), the swept
image of \eqref{eq:rotating-flow-map} is an embedded torus of revolution in
\(H^1\).  The assumption holds whenever $r\ge0$.
More precisely, put
\[
 \widetilde F(u,t)=\bigl(\operatorname{Rot}_{c_\ast t}\beta(u),z(u)\bigr),
 \qquad P_{\rm rot}=2\pi/|c_\ast|.
\]
The rectangular parameter torus
$(\mathbb R/nL\mathbb Z)\times(\mathbb R/P_{\rm rot}\mathbb Z)$
maps $n$-to-one onto the swept torus.
Choose $h\in\mathbb Z$ with $hp\equiv1\pmod n$ and put
\begin{equation}
\label{eq:primitive-torus-return}
 T_{\min}=\frac{2\pi}{n|c_\ast|},\qquad
 \delta=\operatorname{sgn}(c_\ast)hL.
\end{equation}
Then $\widetilde F$ descends to an embedding of
\[
 ((\mathbb R/nL\mathbb Z)\times\mathbb R)/
 \bigl((u,t+T_{\min})\sim(u+\delta,t)\bigr).
\]
This is the primitive mapping-torus parametrization of
Proposition~\ref{thm:cylinder-torus-dichotomy}.
\end{theorem}

\begin{proof}
Write $\rho=|\beta|$. At each $u$ the rotation sweeps out the circle
of radius $\rho(u)$ at height $z(u)$. It therefore suffices first to show
that the meridian $u\mapsto(\rho(u),z(u))$ traces a simple closed curve
in $\rho>0$ over one curvature period $L$.
Put $f(x)=x^2/2+\nu$. By \eqref{eq:torus-condition} and
Theorem~\ref{thm:closing-pairs}, $f(r)<0<f(1)$. In the normalization
of Theorem~\ref{thm:rotating-family}, \eqref{eq:general-rotating-radius} reads
\[
 c_\ast^2\rho^2=K_u^2+f(K)^2=2e_\ast+\nu^2+2c_\ast K.
\]
The first equality shows that $\rho>0$: in the interior of each
half-period $K_u\ne0$, and at either turning value $f(K)\ne0$.
The second equality makes $\rho$ a strictly monotone function of $K$.
Thus each interior radius occurs twice per period, once as $K$ increases
from $r$ to $1$ and once as it decreases back to $r$.

To compare the two heights at that radius, set
\[
 I(x)=\int_r^x\frac{f(y)}{\sqrt{P_{r,\nu}(y)}}\,dy.
\]
By definition $I(r)=0$, and \eqref{eq:rotating-turning} gives
$I(1)=\mathcal B(r,\nu)/2=0$. The function $f$ has exactly one zero
in $(r,1)$: it is increasing for $x\ge0$, and if $r<0$ then
$f(x)\le f(r)<0$ on $[r,0]$. Since the denominator is positive on
$(r,1)$, $I$ first decreases and then increases back to zero.
In particular $I(x)<0$ throughout $(r,1)$.

Let $z_0$ be the height at $K=r$. On the increasing half-period,
$K_u=\sqrt{P_{r,\nu}(K)}$, so integrating
$z_u=-f(K)/c_\ast$ from \eqref{eq:rotating-z-lift} gives
$z=z_0-I(x)/c_\ast$ at $K=x$. At $K=1$ the height is again $z_0$
because $I(1)=0$. On the return half-period $K_u$ changes sign, and
integration from $1$ back to $x$ gives $z=z_0+I(x)/c_\ast$.
Hence the two points at any interior radius have distinct heights.
The branches meet only at the turning points and are regular there
because $z_u=-f(K)/c_\ast\ne0$. Elsewhere $\rho_u\ne0$.
The meridian is therefore smooth, simple and closed. Revolving it about
the axis gives an embedded torus.

It remains to determine how many times $\widetilde F$ covers this image.
If two parameter pairs give the same point, their radii and heights agree.
Since the meridian is embedded over one period, $u_2-u_1=jL$ for some
integer $j$. Advancing by $jL$ rotates $\beta$ through $j\Theta_0$ and
preserves $z$. As $\rho>0$, the points coincide precisely when
\[
 u_2-u_1=jL\pmod{nL},\qquad
 c_\ast(t_2-t_1)+j\Theta_0\in2\pi\mathbb Z.
\]
For each $j=0,\ldots,n-1$ the second condition determines exactly one
time shift modulo $P_{\rm rot}$. Thus every image point has $n$ preimages
on the rectangular torus.

Since $p$ and $n$ are coprime, the angles $j\Theta_0=2\pi jp/n$ modulo
$2\pi$ are precisely the multiples of $2\pi/n$. The shortest positive
time return up to an arclength shift is therefore $T_{\min}=P_{\rm rot}/n$.
With $hp\equiv1\pmod n$, the shift $\delta=\operatorname{sgn}(c_\ast)hL$
rotates the curve through $c_\ast T_{\min}$ modulo $2\pi$, so
$\widetilde F(u,t+T_{\min})=\widetilde F(u+\delta,t)$.
All coincidences above are generated by the two period vectors
$(nL,0)$ and $(-\delta,T_{\min})$. Their fundamental parallelogram has
area $1/n$ times that of the original rectangle. On its quotient the map
is injective. It is an immersion because rotation supplies a nonzero
angular velocity while the meridian has nonzero tangent in $(\rho,z)$.
Compactness of the quotient torus now gives an embedding.
The original $F(s,t)=\widetilde F(s-\nu t,t)$ uses the corresponding
sheared lattice of periods.
The return $t_0=\Theta_0/c_\ast$, $\delta=L$ alone instead gives a
$p$-fold covering when $p>1$, because it omits the shorter rotation returns.

If \(r\ge0\), then \(\mathcal B=0\) forces \(\tfrac12x^2+\nu\) to be negative
somewhere on \([r,1]\), and since it is increasing there,
\(\tfrac12r^2+\nu<0\).
\end{proof}

\begin{corollary}[Embedded tori with prescribed rotation number]
\label{cor:embedded-rotation-numbers}
Let $\Theta^*=\Theta(0,\nu_B(0))$. Every turning angle
$\Theta_0\in[\Theta^*,2\pi)$ occurs on the zero-area graph at some
$0\le r<1$, with $r>0$ if $\Theta_0>\Theta^*$.
If $\Theta_0=2\pi p/n$ with $p/n$ reduced, rotation of the corresponding curve about
the $z$-axis sweeps out an embedded torus by Theorem~\ref{thm:rotating-torus}. In particular, every reduced fraction
$2-2/\sqrt3\le p/n<1$ is realized in this way, including
$(n-1)/n$ for $n\ge7$.
\end{corollary}

\begin{proof}
Along the graph, $r\mapsto\Theta(r,\nu_B(r))$ is continuous on $[0,1]$
and takes the values $\Theta^*$ at $r=0$ and $2\pi$ at $r=1$.
The intermediate value theorem gives the asserted $r$.
Since $r\ge0$, Theorem~\ref{thm:rotating-torus} gives embeddedness
whenever $\Theta_0/(2\pi)$ is rational.

For the explicit bound, take $r=0$ and $\nu=-3/16$. Then
$Q_{0,\nu}(x)=(x+1/2)^2$. Substituting $x=\sin^2 t$ in
\eqref{eq:rotating-turning} and using
$\int_0^{\pi/2}(\sin^2 t+1/2)^{-1}\,dt=\pi/\sqrt3$ gives
\[
 \mathcal B(0,-3/16)=-\frac{\pi}{2\sqrt3},\qquad
 \Theta(0,-3/16)=4\pi(1-1/\sqrt3).
\]
Since $\mathcal B(0,\nu)$ increases with $\nu$, its zero satisfies
$\nu_B(0)>-3/16$. The integrand for $\Theta(0,\nu)$ in
\eqref{eq:rotating-turning} decreases with $\nu$ on $0<x<1$.
Hence $\Theta^*/(2\pi)<2-2/\sqrt3<6/7$.
\end{proof}

\begin{remark}[Monotonicity on the flat zero-area graph]
\label{rem:flat-rotation-monotonicity}
Set $R(r)=\Theta(r,\nu_B(r))/(2\pi)$. Corollary~\ref{cor:embedded-rotation-numbers}
uses continuity and the endpoint values of $R$ to obtain a normalized pair
with $0\le r<1$ for each prescribed rotation number in $[R(0),1)$.
Uniqueness on this branch follows from the interval estimates below.

The following bounds are computer-assisted.  We bracket $\nu_B$ on each
parameter cell by opposite signs of $\mathcal B$, differentiate \eqref{eq:rotating-turning-regular} along
$\mathcal B=0$ using $\nu_B'=-\mathcal B_r/\mathcal B_\nu$, and enclose
the integrals in \eqref{eq:rotating-turning-regular} and
\eqref{eq:rotating-area-regular} by outward-rounded interval sums on
$1024$ equal angular cells.  This gives
\[
\begin{aligned}
 \nu_B(r)+r^2/2 &> 0.00072 &&(-1\le r\le-1/2),\\
 R'(r) &> 0.00022 &&(-1/2\le r\le7/8),\\
 R''(r) &< -0.0114 &&(7/8\le r\le1).
\end{aligned}
\]
At $r=1$, $x_r(\phi)=1$, $Q=4$, and $\partial_r x_r=(1-\cos\phi)/2$
has mean $1/2$. Differentiating the area and angle integrals there gives
$\nu_B'(1)=-1/2$ and $R'(1)=0$. The last bound thus gives $R'(r)>0$
for $7/8\le r<1$. Together with the bound on $R'$ this proves strict
increase on $[-1/2,1)$ and hence the stated uniqueness. By the first bound, the entire thick segment of
Figure~\ref{fig:flat-closing-map} lies in $-1/2<r<1$.

Numerically, $\nu_B(0)\simeq-0.13744$ and $R(0)\simeq0.78039$, below $4/5$.
The $1/2$ pair at $r\simeq-0.2615$ also satisfies \eqref{eq:torus-condition};
Figure~\ref{fig:flat-closing-map} shows the crossing of its boundary near $r=-0.49$.
\end{remark}

\begin{remark}[Dependence of the swept surface on $c_0$]
\label{rem:torus-reeb-component}
The embedded torus is the image of the rotation $F$ in
\eqref{eq:rotating-flow-map}, whose velocity has Reeb component $(2K+c_0)\xi$.
To obtain the representative with Reeb component $2K\xi$ in
\eqref{eq:sasakian-curve-flow}, translate $F$ vertically by $-c_0t$:
\[
 F_0(s,t)=\bigl(\operatorname{Rot}_{c_\ast t}\beta(s-\nu t),z(s-\nu t)-c_0t\bigr).
\]
If $c_0\ne0$, each full rotation has a nonzero vertical displacement,
and the time parameter is noncompact.  This does not imply that the swept
image is an embedded cylinder: the map can have self-intersections.
The maps $F$ and $F_0$ represent the same Hamiltonian trajectory on
$\mathcal Q$, and they coincide when $c_0=0$.
\end{remark}

\subsection{Closed relative equilibria in the special linear group}
\label{sec:split-relative-equilibrium}

We now apply Corollary~\ref{thm:sasaki-projection} and
Theorem~\ref{thm:unified-holonomy} to the Riemannian split member $\lambda=-2$.
Its base is the hyperbolic plane $\mathbb H^2(-4)$.

Use the hyperboloid model, writing a point as $p\in\mathbb R^{2,1}$ with
$\langle p,p\rangle=-1$.  Since the metric of $\mathbb H^2(-4)$ is one
quarter of the unit-curvature metric, take a Lorentz-orthonormal frame
$(p,T,N)$.  The unit tangent and normal for the curvature $-4$ metric are
$2T$ and $2N$, respectively, and the arclength Frenet equations are
\begin{equation}
\label{eq:hyperbolic-frenet}
 p_s=2T,\qquad T_s=2p+kN,\qquad N_s=-kT.
\end{equation}
Equivalently, for $F=(p,T,N)\in SO^+(2,1)$,
\begin{equation}
\label{eq:hyperbolic-frame-matrix}
 F_s=F A(k),\qquad
 A(k)=
 \begin{pmatrix}
 0&2&0\\ 2&0&-k\\ 0&k&0
 \end{pmatrix}.
\end{equation}

The monodromy criterion \eqref{eq:space-form-base-closing} and
Proposition~\ref{prop:space-form-noether-reconstruction} apply directly in
this model.  We therefore restrict attention to the explicit evaluation of
the frame monodromy.

\subsubsection*{The branch \(C=0\): closure after one curvature period}
For the symmetric branch \(C=0\), write
\begin{equation}
\label{eq:split-stationary-ode}
 K''+\frac12K^3+\alpha K=0.
\end{equation}
Since \(\kappa=2\lambda=-4\), comparison with
\eqref{eq:space-form-profile-again} gives the traveling-wave speed
\[
 v_{\mathrm{hyp}}=\alpha+4.
\]

\begin{proposition}[The boost increment over one curvature period]
\label{prop:split-boost-function}
Consider the Jacobi curvature profile
\begin{equation}
\label{eq:split-cn-family}
 K_{m,\omega}(s)=2\sqrt m\,\omega\,
 \cn(\omega s;m),\qquad L=\frac{4\mathbb K(m)}{\omega},
\end{equation}
which is the \(C=0\) specialization of
Proposition~\ref{prop:general-cn-wave} for \(\lambda=-2\).  It solves
\eqref{eq:split-stationary-ode} with
$\alpha=(1-2m)\omega^2$.  Along this solution the ambient vector
\begin{equation}
\label{eq:split-conserved-vector}
 V=Kp-\frac12K_sT+
 \left(\frac14K^2+\frac{\alpha+4}{2}\right)N
\end{equation}
is constant. Its squared Minkowski norm is
\begin{equation}
\label{eq:split-D}
 \mathfrak d_V=\langle V,V\rangle
 =-4m\omega^2+\left(\frac12\omega^2+2\right)^2.
\end{equation}
When $\mathfrak d_V>0$, let $\Lambda_\delta$ denote the Lorentz boost fixing
$V$ with signed hyperbolic angle $\delta$, oriented as in the proof.
At fixed time, one curvature period takes the frame $F(s)$ to
$F(s+L)=\Lambda_\delta F(s)$ and hence the point $p(s)$ to $p(s+L)=\Lambda_\delta p(s)$.
The angle accumulated over this period is
\begin{equation}
\label{eq:split-boost-increment}
 \Delta_m(\omega)=
 \int_0^L
 \frac{2\sqrt{\mathfrak d_V}\left(\frac14K^2+\frac{\alpha+4}{2}\right)}
 {\mathfrak d_V+K^2}\,ds,
\end{equation}
so $\delta=\Delta_m(\omega)$. Since $\Lambda_\delta=I$ exactly when $\delta=0$,
the frame returns after one curvature period if and only if
$\Delta_m(\omega)=0$.

For the time evolution, define $A_V\in\mathfrak{so}(2,1)$ by
$A_V X=V\boxtimes X$, where $\boxtimes$ is the Lorentzian cross product
with $\langle a\boxtimes b,c\rangle=\det(a,b,c)$.
After the arclength shift is removed, the projected curve moves by
$\exp(4tA_V)$, a hyperbolic one-parameter group when $V$ is spacelike.

\begin{proof}
Differentiating \eqref{eq:split-conserved-vector}, using
\eqref{eq:hyperbolic-frenet} and \eqref{eq:split-stationary-ode}, gives
$V_s=0$.  Evaluating its norm at a maximum of $K$ gives
\eqref{eq:split-D}.  The stabilizer of a spacelike vector in
$SO^+(2,1)$ is the one-dimensional boost group.  Projecting $p$ to
$V^\perp$, put
\[
 p_\perp=p+\frac K{\mathfrak d_V} V,
 \qquad
 -\langle p_\perp,p_\perp\rangle=\frac{\mathfrak d_V+K^2}{\mathfrak d_V}.
\]
The oriented area form on $V^\perp$ gives
\[
 \frac{\det(V,p_\perp,(p_\perp)_s)/\sqrt{\mathfrak d_V}}
      {-\langle p_\perp,p_\perp\rangle}
 =\frac{2\sqrt{\mathfrak d_V}\left(\frac14K^2+\frac{\alpha+4}{2}\right)}{\mathfrak d_V+K^2}.
\]
This is the derivative of the Lorentzian polar angle and proves
\eqref{eq:split-boost-increment} upon integration over $[0,L]$.
Finally, the projected mKdV velocity, expressed in
the ambient hyperboloid coordinates, satisfies
\[
 p_t+v_{\mathrm{hyp}}p_s=(K^2+2v_{\mathrm{hyp}})T+2K_sN
             =4A_Vp,
\]
which gives the asserted time evolution.
\end{proof}
\end{proposition}

\begin{theorem}[Analytic split-real closing theorem]
\label{thm:analytic-split-closure}
Let $m_*$ be the elliptic parameter in Theorem~\ref{thm:closed-elliptic}.  For every
\[
 m\in(m_*,1)
\]
there is a unique $\omega_m>0$ for which the Jacobi curvature
$K_{m,\omega_m}$ in \eqref{eq:split-cn-family} has identity frame
monodromy, and $\omega_m>2$. The function $m\mapsto\omega_m$ is real analytic on $(m_*,1)$.
Its projected curve in $\mathbb H^2(-4)$ is smoothly closed,
has zero signed area, and is the union of an arc and its image under a
hyperbolic reflection.  Consequently its horizontal lift to
$\widetilde{SL}(2,\mathbb R)$ is a closed Legendrian relative equilibrium.
For $m\in(0,m_*]$ no $\omega>0$ gives identity frame monodromy.
\end{theorem}

\begin{proof}
Put $x=\omega^2$ and define
\[
 \ell_m(x)=\frac{16mx}{(x+4)^2}.
\]
For $x\geq4$ one has $0<\ell_m(x)\leq m<1$.  Let
\[
 \boldsymbol\Pi(\ell_\Pi;m)
 =\int_0^{\mathbb K(m)}
   \frac{du}{1-\ell_\Pi\sn^2(u;m)}
\]
be the complete elliptic integral of the third kind.  Along
\eqref{eq:split-cn-family},
\begin{align*}
 \frac14K^2+\frac{\alpha+4}{2}
   &=\frac{x+4}{2}-mx\sn^2(u;m),\\
 \mathfrak d_V+K^2
   &=\frac{(x+4)^2}{4}-4mx\sn^2(u;m),
 \qquad u=\omega s.
\end{align*}
The elementary decomposition
\[
 \frac{\frac{x+4}{2}-mx\sn^2u}
      {\frac{(x+4)^2}{4}-4mx\sn^2u}
 =\frac14+
 \frac{(x+4)(4-x)}
 {16\left(\frac{(x+4)^2}{4}-4mx\sn^2u\right)}
\]
and integration over two periods of $\sn^2$ turn
\eqref{eq:split-boost-increment} into the exact formula
\begin{equation}
\label{eq:split-Pi-closing}
 \Delta_m(\sqrt x)
 =\frac{2\sqrt{\mathfrak d_V(x)}}{\sqrt x}\,F_m(x),
 \qquad
 F_m(x)=\mathbb K(m)
 +\frac{4-x}{4+x}\,
   \boldsymbol\Pi\!\left(\ell_m(x);m\right).
\end{equation}
Here
\begin{equation}
\label{eq:split-D-positive}
 \mathfrak d_V(x)
 =\frac14\bigl(x-(8m-4)\bigr)^2+16m(1-m)>0.
\end{equation}
Thus the prefactor in \eqref{eq:split-Pi-closing} never vanishes.

For $0<x\le4$ the coefficient $(4-x)/(4+x)$ is nonnegative and
$\boldsymbol\Pi>0$, so $F_m(x)\ge\mathbb K(m)>0$; every zero of $F_m$
therefore satisfies $x>4$, that is, $\omega>2$.  The function $F_m$ is
continuous on $[4,\infty)$ and $F_m(4)=\mathbb K(m)>0$.  To determine its sign at infinity, use
\[
 \int_0^{\mathbb K(m)}\sn^2(u;m)\,du
 =\frac{\mathbb K(m)-\mathbb E(m)}{m}.
\]
Uniform expansion of the defining integral for $\boldsymbol\Pi$ at $\ell_\Pi=0$
gives
\[
 \boldsymbol\Pi(\ell_\Pi;m)
 =\mathbb K(m)
  +\ell_\Pi\,\frac{\mathbb K(m)-\mathbb E(m)}m+O(\ell_\Pi^2).
\]
Since
\[
 \ell_m(x)=\frac{16m}{x}+O(x^{-2}),
 \qquad
 \frac{4-x}{4+x}=-1+\frac8x+O(x^{-2}),
\]
we obtain
\begin{equation}
\label{eq:split-F-asymptotic}
 F_m(x)
 =\frac{8\bigl(2\mathbb E(m)-\mathbb K(m)\bigr)}x
  +O(x^{-2}).
\end{equation}
The ratio $\mathbb E/\mathbb K$ is strictly decreasing, as proved in
Theorem~\ref{thm:closed-elliptic}.  Therefore
$2\mathbb E(m)-\mathbb K(m)<0$ whenever $m>m_*$.  It follows that
$F_m(x)<0$ for all sufficiently large $x$.  The intermediate value theorem
now supplies $x_m>4$ with $F_m(x_m)=0$.  Setting
$\omega_m=\sqrt{x_m}$ and using
Proposition~\ref{prop:split-boost-function} gives identity monodromy and
hence a smoothly closed projected curve.

To prove uniqueness and analytic dependence, set $t=(x-4)/(x+4)$,
$\ell=m(1-t^2)$, and write $K=\mathbb K(m)$, $E=\mathbb E(m)$ and
$G(t)=K-t\boldsymbol\Pi(\ell;m)$, so that $F_m(x)=G(t)$ and
$\ell=\ell_m(x)$. For $0<t<1$, the characteristic
derivative of the complete elliptic integral is
\[
 2(m-\ell)(\ell-1)\ell\,\partial_\ell\boldsymbol\Pi(\ell;m)
 =\ell E+(m-\ell)K+(\ell^2-m)\boldsymbol\Pi(\ell;m).
\]
This identity follows by differentiating the defining integral and
integrating the derivative of
$\sin\theta\cos\theta\sqrt{1-m\sin^2\theta}/(1-\ell\sin^2\theta)$
over $[0,\pi/2]$, where its boundary values vanish.
Using $\ell'=-2mt$ gives
\[
 (1-\ell)G'(t)+\frac{t}{1-t^2}G(t)=\frac{t}{1+t}K-E .
\]
Put
\[
 \mu(t)=\exp\int_0^t\frac{y\,dy}{(1-y^2)(1-\ell(y))},
 \qquad\text{so that}\qquad
 (\mu G)'=\frac{\mu}{1-\ell}\Bigl(\frac{t}{1+t}K-E\Bigr).
\]
Here $1-\ell\ge1-m>0$.  Since $E/K<1/2$, the right side changes sign
exactly once, from negative to positive, so $\mu G$ first decreases
strictly and then increases strictly. Now $G(0)=K>0$, while
\eqref{eq:split-F-asymptotic} gives $G(t)<0$ for $t$ sufficiently close
to $1$. Hence the sole interior zero lies on the decreasing interval and
is simple. The integral is jointly real analytic in $(m,t)$ there.
The analytic implicit function theorem, followed by
$\omega_m=2\sqrt{(1+t)/(1-t)}$, gives the assertion.

Now let $0<m\le m_*$, so that $E/K\ge1/2$.  Since $t/(1+t)<1/2$, the
right side above is negative on $(0,1)$, and $\mu G$ decreases strictly
there.  As $t\to1$, $1-\ell\to1$, so $\mu(t)=O\bigl((1-t)^{-1/2}\bigr)$,
while \eqref{eq:split-F-asymptotic} and $x=4(1+t)/(1-t)$ give
$G(t)=O(1-t)$.  Hence $\mu G\to0$, and strict decrease gives $\mu G>0$ on
$(0,1)$.  Thus $F_m$ has no zero with $x>4$, and none with $0<x\le4$ by
the first step, so no $\omega>0$ gives identity monodromy.

It remains to prove zero Reeb holonomy.
Let $Y_s=Y A(K)$, $Y(0)=I$, and put
\[
 h=\frac L2,\qquad H=Y(h),\qquad
 S=\operatorname{diag}(1,1,-1).
\]
The half-period identity $K(s+h)=-K(s)$ and
$A(-K)=SA(K)S$ imply
\[
 \mathsf M=HSHS=(HS)^2.
\]
At the zero found above, $\mathsf M=I$.  If $F(0)$ is the initial Lorentz frame,
then
\[
 \mathscr S=F(0)HSF(0)^{-1}
\]
is an orientation-reversing involutive isometry of $\mathbb H^2(-4)$,
namely a hyperbolic reflection, and
\[
 p(s+h)=\mathscr S p(s).
\]
Thus the second half of the closed curve is the reflected image of the
first half with the same parameter orientation.  To see the area
cancellation directly, join the endpoints of the first arc by the
$\mathscr S$-invariant geodesic segment and let $\mathscr D_+$ be the resulting signed
two-chain.  The full curve bounds $\mathscr D_++\mathscr S(\mathscr D_+)$, while $\mathscr S$ reverses the
hyperbolic area form.  Consequently
\[
 \operatorname{Area}_{\rm sgn}\bigl(\mathscr D_++\mathscr S(\mathscr D_+)\bigr)=0.
\]
By Theorem~\ref{thm:unified-holonomy}, the horizontal lift closes.
Its relative-equilibrium property follows from the general
Noether--Killing reconstruction, Proposition~\ref{prop:space-form-noether-reconstruction}.
\end{proof}

\begin{figure}[!htbp]
 \centering
 \includegraphics[width=\linewidth]{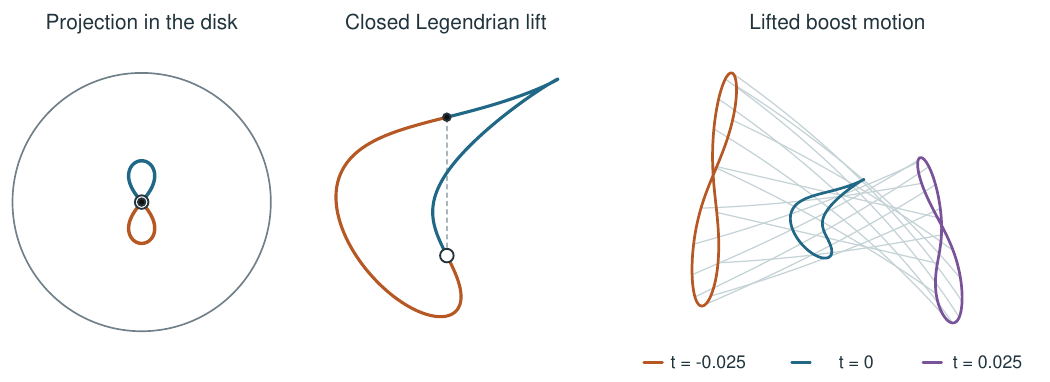}
 \caption{The closed boost family of Theorem~\ref{thm:analytic-split-closure}
 at $m=0.85$ and the frequency in \eqref{eq:split-numerical-root}.
 Left: the two reflected halves of the curve in the Poincar\'e disk.
 The open circle and solid dot coincide at the crossing. Middle: the
 corresponding points have distinct heights on the closed horizontal lift
 in disk and fiber coordinates. Right: three positions of this lift under
 the lifted boost in $\widetilde{SL}(2,\mathbb R)$, with the lift chosen so that
 $\eta(\widehat Y)=2K$ along the initial curve. Thin lines trace point
 orbits. Both spatial views magnify the fiber coordinate by a factor of $2.8$.}
 \label{fig:split-real-geometry}
\end{figure}

\begin{remark}
For $m=0.85$, numerical solution of \eqref{eq:split-Pi-closing} gives
\begin{equation}
\label{eq:split-numerical-root}
 \omega=5.64725040011834\ldots,\qquad
 L=1.69216260449546\ldots.
\end{equation}
These values are used for Figure~\ref{fig:split-real-geometry}.
\end{remark}

\begin{proposition}[Parameter topology of a swept relative equilibrium]
\label{thm:cylinder-torus-dichotomy}
Let $\gamma_0:\mathbb R/\ell\mathbb Z\to G_\lambda$ be a closed
Legendrian immersion and write $\gamma(s,t)=\Phi_t\gamma_0(s-vt)$.
If the orbit map $\widetilde\gamma(u,t)=\Phi_t\gamma_0(u)$ has rank two,
it is an immersion of the cylinder
$(\mathbb R/\ell\mathbb Z)\times\mathbb R$.
If there are constants $t_0>0$ and $\delta$ such that
\begin{equation}
\label{eq:time-shift-identification}
 \Phi_{t_0}\gamma_0(u)=\gamma_0(u+\delta)\quad\text{for every }u,
\end{equation}
the orbit map descends to
\[
 ((\mathbb R/\ell\mathbb Z)\times\mathbb R)/
 \bigl((u,t+t_0)\sim(u+\delta,t)\bigr),
\]
a mapping torus of a circle translation, hence a torus.
The cylindrical map is an embedding if it is injective and proper;
the descended torus map is an embedding if it is injective.
A time return modulo reparametrization need not be a literal period
with the original parameter $s$ fixed.

Let $\lambda\le0$, and suppose that the projected generator is a nonzero
Killing field $Y$ of $\Sigma_{2\lambda}$ that is not elliptic: a
translation for $\lambda=0$, or a parabolic or hyperbolic field for
$\lambda<0$.  Then $\gamma_0$ has no nonzero time return, even setwise,
and the time parameter stays cylindrical for every choice of $c_0$.
This applies to the figure-eight of Theorem~\ref{thm:closed-elliptic} and
to the boost family of Theorem~\ref{thm:analytic-split-closure}.  No
embeddedness assertion is made for these swept maps.
\end{proposition}

\begin{proof}
The substitution $u=s-vt$ is a diffeomorphism of the cylinder.
Equation~\eqref{eq:time-shift-identification} gives
$\widetilde\gamma(u,t+t_0)=\widetilde\gamma(u+\delta,t)$.
The coordinate $u+\delta t/t_0$ is invariant under
$(u,t)\mapsto(u-\delta,t+t_0)$, identifying the quotient with a torus.
The embedding criteria follow from proper injective immersion and
compactness, respectively.

A time return, even setwise, would give $t_0\ne0$ with
$\exp(t_0Y)\beta(S^1)=\beta(S^1)$ for the projected curve $\beta$; a
constant Reeb translation does not change the projection.  For a nonzero
translation of the plane, and for a nonzero parabolic or hyperbolic
isometry of the hyperbolic plane, the iterates of every point leave every
compact set.  Hence $\exp(t_0Y)$ cannot preserve the nonempty compact set
$\beta(S^1)$.  The figure-eight is generated by a translation, since $C=0$
(Remark~\ref{rem:noether-criterion}).  For the split family,
$\mathfrak d_V>0$ in \eqref{eq:split-D-positive}, so the projected
generator is a nontrivial boost (Proposition~\ref{prop:split-boost-function}).
\end{proof}

\subsubsection*{The branch \(C\ne0\): elliptic rotating relative equilibria}

We keep $\lambda=-2$, $\kappa=-4$, and the normalization of
Section~\ref{sec:traveling-waves}: $K(s)=b\,w(bs)$, where $w$ oscillates
between $r$ and $1$ with $w(0)=1$, and $(r,\nu)\in\mathcal{A}$.  Unlike the
flat case, the hyperbolic base breaks the scaling symmetry, so the scale
enters through
\begin{equation}
\label{eq:split-eps}
 \varepsilon=\frac4{b^2}>0,\qquad v=b^2(\nu+\varepsilon)
\end{equation}
by \eqref{eq:denormalization}.  Along the reconstructed curve,
Proposition~\ref{prop:space-form-noether-reconstruction} gives
$\psi=C+4K=b^3(c_\ast+\varepsilon w)$ and
\begin{equation}
\label{eq:split-rot-invariants}
 |Y|^2=b^4\bigl(w_u^2+(\tfrac12w^2+\nu+\varepsilon)^2\bigr),\qquad
 \mathcal I=b^6\varpi_\varepsilon^2,\qquad
 \varpi_\varepsilon^2:=c_\ast^2-\varepsilon\bigl(2e_\ast+(\nu+\varepsilon)^2\bigr).
\end{equation}
Thus the generator is elliptic exactly when $\varpi_\varepsilon^2>0$.  In
that case we write $\varpi_\varepsilon>0$ for the positive root and, for
$c_\ast>0$, deform \eqref{eq:rotating-area-function} to
\begin{equation}
\label{eq:split-area-function}
 \mathcal B_\varepsilon(r,\nu)
 =\int_0^{L_0}\frac{2c_\ast(\frac12w^2+\nu+\varepsilon)}
 {c_\ast+\varepsilon w+\varpi_\varepsilon}\,du
 =4\int_0^\pi
 \frac{2c_\ast\bigl(\frac12x_r(\phi)^2+\nu+\varepsilon\bigr)}
 {\bigl(c_\ast+\varepsilon x_r(\phi)+\varpi_\varepsilon\bigr)
  \sqrt{Q_{r,\nu}(x_r(\phi))}}\,d\phi .
\end{equation}
Since $\varpi_0=c_\ast$, one has $\mathcal B_0=\mathcal B$.

\begin{lemma}[Rotating closing data on $\mathbb H^2(-4)$]
\label{lem:split-rotating-data}
Let $(r,\nu)\in\mathcal{A}$ with $r<1$, $c_\ast>0$, $c_\ast+\varepsilon r>0$, and
$\varpi_\varepsilon^2>0$.  Then the Killing field $Y$ of
Proposition~\ref{prop:space-form-noether-reconstruction} is the rotation
$Y=\omega_0\partial_\varphi$, $\omega_0=b^3\varpi_\varepsilon$, about a point
$o$, where $(\rho,\varphi)$ are geodesic polar coordinates at $o$ and the
metric is $d\rho^2+\frac14\sinh^2(2\rho)\,d\varphi^2$.  Along the curve
\begin{equation}
\label{eq:split-rot-radius}
 \cosh2\rho=\frac{c_\ast+\varepsilon w}{\varpi_\varepsilon}.
\end{equation}
The monodromy $\mathcal M_K$ is the rotation about $o$ through an angle
$\varphi_L$ with
\begin{equation}
\label{eq:split-rot-angle}
 \varphi_L\equiv\Theta(r,\nu)-\frac{\varepsilon}{2c_\ast}\,\mathcal B_\varepsilon(r,\nu)
 \pmod{2\pi},
\end{equation}
and the signed area contributed by one curvature period is
$\mathcal B_\varepsilon(r,\nu)/(2b^2c_\ast)$.  Consequently, if
\begin{equation}
\label{eq:split-rot-closing}
 \Theta(r,\nu)=\frac{2\pi p}{n},\qquad \mathcal B_\varepsilon(r,\nu)=0,
\end{equation}
with $p/n$ reduced, then the projected curve closes after $n$ curvature
periods with zero signed area, and its horizontal lift to
$\widetilde{SL}(2,\mathbb R)$ is a closed Legendrian relative equilibrium
generated by the lift of the elliptic field $Y$.
\end{lemma}

\begin{proof}
The rotation field $\partial_\varphi$ has Killing data $(0,1)$ at $o$ and
satisfies $\nabla\partial_\varphi=\cosh(2\rho)\,J$.  By
\eqref{eq:killing-casimir}, $\psi^2=\mathcal I+4|Y|^2\ge\mathcal I>0$, so
$\psi$ has constant sign on $\mathbb H^2(-4)$; it is positive on the curve
because $c_\ast+\varepsilon w\ge c_\ast+\varepsilon r>0$.  Since an elliptic
Killing field vanishes at one point $o$ and $\psi(o)^2=\mathcal I$, we get
$Y=\omega_0\partial_\varphi$ with $\omega_0=b^3\varpi_\varepsilon$ and
$\psi=\omega_0\cosh2\rho$, which is \eqref{eq:split-rot-radius}.

The $1$-form $\alpha=\frac14(\cosh2\rho-1)\,d\varphi$ is smooth on
$\mathbb H^2(-4)$, rotation invariant, and satisfies
$d\alpha=\frac12\sinh(2\rho)\,d\rho\wedge d\varphi=dA_{-4}$.  Where
$\rho>0$, $\varphi_s=\langle T,\partial_\varphi\rangle/|\partial_\varphi|^2
=\omega_0(\frac12K^2+v)/|Y|^2$, and $\psi^2-\omega_0^2=4|Y|^2$ gives
\[
 \alpha(\beta_s)=\frac{(\psi-\omega_0)(\frac12K^2+v)}{4|Y|^2}
 =\frac{\frac12K^2+v}{\psi+\omega_0}.
\]
Both sides are continuous, so this holds everywhere.  Integrating over one
period and substituting \eqref{eq:split-eps} gives the area
$\mathcal B_\varepsilon/(2b^2c_\ast)$.

Since $\mathcal M_K$ maps $Y$ along $\beta$ to itself, it commutes with
the one-parameter group of $Y$ and is a rotation about $o$.  Where
$\rho>0$, let $\theta$ be the angle from $\partial_\rho$ to $T$.  The
connection form of the polar frame is $\cosh(2\rho)\,d\varphi$, hence
$K=\theta_s+\cosh(2\rho)\varphi_s$.  The rotation preserves the polar frame,
so $\theta(L)\equiv\theta(0)$, and
\[
 \int_0^LK\,ds\equiv\int_0^L\cosh(2\rho)\,\varphi_s\,ds
 =\varphi_L+4\int_0^L\alpha(\beta_s)\,ds\pmod{2\pi}.
\]
Since $\int_0^LK\,ds=\Theta(r,\nu)$ and $4b^{-2}=\varepsilon$, this is
\eqref{eq:split-rot-angle}.  The curve passes through $o$ only when
$\frac12K^2+v$ vanishes at a turning value, a codimension-one condition, and
both sides of \eqref{eq:split-rot-angle} are continuous modulo $2\pi$; hence
the identity holds without this restriction.

Under \eqref{eq:split-rot-closing}, $\mathcal M_K$ has order $n$, so the
base curve closes after $n$ periods by \eqref{eq:space-form-base-closing}.
Each period is the image of the first under a power of $\mathcal M_K$,
which preserves $\alpha$, so the total signed area is
$n\mathcal B_\varepsilon/(2b^2c_\ast)=0$.  Theorem~\ref{thm:unified-holonomy}
closes the lift, and Proposition~\ref{prop:space-form-noether-reconstruction}
with Theorem~\ref{thm:relative-equilibria} gives the relative-equilibrium
property.
\end{proof}

\begin{figure}[!t]
 \centering
 \includegraphics[width=.97\linewidth]{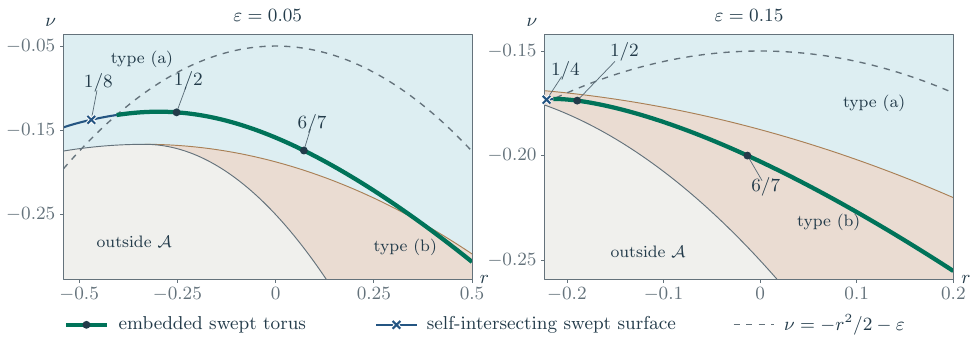}
 \caption{Elliptic rotations in $\widetilde{SL}(2,\mathbb R)$. Each solid curve
 is a computed segment of $\mathcal B_\varepsilon=0$. Each label $p/n$
 gives the rotation per curvature period in turns. Successive arcs of one
 curvature period differ by a rotation of $2\pi p/n$, and the lift closes
 after $n$ periods. For the lifted rotation at fixed fiber height, dots below the
 dashed parabola give embedded swept tori and crosses above it give self-intersections.
 Shading distinguishes the curvature profiles. Figure~\ref{fig:split-rotating-torus}
 shows both $6/7$ tori; the one at $\varepsilon=0.05$ is of type (a) and the one at $\varepsilon=0.15$ of type (b).}
 \label{fig:split-parameter-map}
\end{figure}

At a fixed scale on the elliptic branch, we likewise call $(r,\nu)$ a
\emph{closing pair} when it satisfies \eqref{eq:split-rot-closing} for a
reduced fraction $p/n$. Figure~\ref{fig:split-parameter-map} plots these
conditions at two scales.

Thus the first closing equation is the flat equation unchanged; only the
area condition is deformed.  The deformation is regular in the flat limit.

\begin{lemma}[Flat limit of the area function]
\label{lem:split-flat-limit}
Let $R\subset\mathcal{A}$ be compact with $r<1$ and $c_\ast>0$ on $R$.  There
is $\varepsilon_R>0$ such that for $0<\varepsilon\le\varepsilon_R$ the
hypotheses of Lemma~\ref{lem:split-rotating-data} hold on $R$, and
\[
 \mathcal B_\varepsilon\longrightarrow\mathcal B
 \qquad(\varepsilon\downarrow0)
\]
uniformly on $R$.
\end{lemma}

\begin{proof}
On $R$, $c_\ast$ and $Q_{r,\nu}(x_r(\phi))$ are bounded below by positive
constants and $e_\ast,\nu$ are bounded.  Hence $\varpi_\varepsilon\to c_\ast$
and $c_\ast+\varepsilon x_r+\varpi_\varepsilon\to2c_\ast$ uniformly, and the
integrand of $\mathcal B_\varepsilon$ converges uniformly to that of
\eqref{eq:rotating-area-regular}.
\end{proof}

\begin{theorem}[Closed rotating relative equilibria in $\widetilde{SL}(2,\mathbb R)$]
\label{thm:split-rotating}
For every $\Theta_0\in(0,2\pi)$ there is $\varepsilon_0>0$ such that for
every $\varepsilon\in(0,\varepsilon_0)$ there is $(r,\nu)\in\mathcal{A}$ with
$-1<r<1$, $c_\ast>0$, $\varpi_\varepsilon^2>0$,
\[
 \Theta(r,\nu)=\Theta_0,\qquad \mathcal B_\varepsilon(r,\nu)=0.
\]
Consequently, for every reduced fraction $p/n\in(0,1)$ there is $b_0>0$
such that for every $b>b_0$ the curvature $K=b\,w(bs)$ of the corresponding
profile reconstructs a closed Legendrian curve in $\widetilde{SL}(2,\mathbb R)$,
closing after $n$ curvature periods, which the Hamiltonian flow of
horizontal length moves as a relative equilibrium generated by an elliptic
Killing field.
If $\Theta_0\in(\Theta^*,2\pi)$, with $\Theta^*$ as in
Corollary~\ref{cor:embedded-rotation-numbers}, the pair can be chosen with $0<r<1$.
\end{theorem}

\begin{proof}
By the endpoint values of $\vartheta$, choose $-1<r_-<r_+<1$ such that
$\vartheta(r_-)<\Theta_0<\vartheta(r_+)$. If $\Theta_0>\Theta^*$ take
$r_-=0$, since $\vartheta(0)=\Theta^*$. On the compact graph segment over
$I=[r_-,r_+]$, Lemma~\ref{lem:area-graph} gives $c_\ast>0$ and
$\partial_\nu\mathcal B>0$. Both have positive lower bounds there.

In \eqref{eq:split-area-function} the radicand defining
$\varpi_\varepsilon$ equals $c_\ast^2$ at $\varepsilon=0$.
Its positive square root, the denominator of the integrand, and
$Q_{r,\nu}$ therefore stay positive on a fixed neighborhood of the
segment for sufficiently small $|\varepsilon|$. The integral is jointly
real analytic in $(r,\nu,\varepsilon)$ there, and
$\mathcal B_0=\mathcal B$. The analytic implicit function theorem
continues each point of the graph to a local solution
$\nu=\nu_{B,\varepsilon}(r)$. By compactness we may take a common
$\varepsilon$-interval for $r\in I$. The local solutions agree on
overlaps: in a sufficiently narrow tube about the flat graph,
$\partial_\nu\mathcal B_\varepsilon$ stays positive. This gives a
single perturbed graph on $I$, converging uniformly to $\nu_B$.
By Lemma~\ref{lem:split-flat-limit} its points satisfy the admissibility
conditions of Lemma~\ref{lem:split-rotating-data} for small positive
$\varepsilon$.

The strict inequalities at $r_-,r_+$ persist for
$r\mapsto\Theta(r,\nu_{B,\varepsilon}(r))$. The intermediate value
theorem gives the required angle in $(r_-,r_+)$, hence at positive $r$
if $r_-=0$. For $\Theta_0=2\pi p/n$, Lemma~\ref{lem:split-rotating-data}
closes the lift after $n$ curvature periods. Substituting
$\varepsilon=4/b^2$ gives the stated range of scales.
\end{proof}

\begin{theorem}[Embedded swept tori in $\widetilde{SL}(2,\mathbb R)$]
\label{thm:split-torus}
Suppose the hypotheses of Lemma~\ref{lem:split-rotating-data} hold.
Assume \eqref{eq:split-rot-closing} and
\begin{equation}
\label{eq:split-torus-condition}
 \tfrac12r^2+\nu+\varepsilon<0.
\end{equation}
Then the lifted rotation fixing the fiber over the elliptic center pointwise
sweeps out an embedded torus in
$\widetilde{SL}(2,\mathbb R)$.  Condition \eqref{eq:split-torus-condition}
holds whenever $r\ge0$; in particular it holds for the pairs of
Theorem~\ref{thm:split-rotating} with $\Theta_0\in(\Theta^*,2\pi)$.
\end{theorem}

\begin{proof}
The Reeb fibration over the contractible base is a trivial $\mathbb
R$-bundle with connection $\eta$ and curvature $d\eta=2\pi_\lambda^*dA$, and
the Sasakian metric is $\pi_\lambda^*g_{\mathbb H^2}+\eta\otimes\eta$.  After
changing the fiber coordinate by $z\mapsto z+f$ we may therefore use coordinates
$(\rho,\varphi,z)$ with
\[
 \eta=dz+\tfrac12(\cosh2\rho-1)\,d\varphi=dz+2\alpha.
\]
The maps $(\rho,\varphi,z)\mapsto(\rho,\varphi+\omega_0t,z)$ are Sasakian
automorphisms lifting the flow of $Y$, with potential
$f_Y=\frac12\omega_0(\cosh2\rho-1)$.  Along the curve
$f_Y\circ\beta=\frac12(\psi-\omega_0)=2K+\frac12(C-\omega_0)$, so by
Theorem~\ref{thm:relative-equilibria} these rotations realize the relative
equilibrium for $c_0=\frac12(C-\omega_0)$.  The swept image is
therefore the surface of revolution of the meridian $(\rho(s),z(s))$, with
\[
 z_s=-2\alpha(\beta_s)=-\frac{2(\frac12K^2+v)}{\psi+\omega_0}.
\]
By \eqref{eq:split-rot-radius}, $\cosh2\rho$ is a strictly increasing
affine function of $w$.  At the turning values, $\frac12K^2+v=b^2(\frac12w^2+\nu+\varepsilon)$
is negative at $w=r$ by \eqref{eq:split-torus-condition} and positive at
$w=1$, because $\mathcal B_\varepsilon=0$ forces a sign change.  Hence
$|Y|>0$ and $\rho>0$.  The function $\frac12x^2+\nu+\varepsilon$ changes
sign exactly once on $[r,1]$, as in the proof of
Theorem~\ref{thm:rotating-torus}, and the weight
$(c_\ast+\varepsilon x+\varpi_\varepsilon)^{-1}$ is positive.  Therefore
\[
 I_\varepsilon(x)=\int_r^x
 \frac{\frac12u^2+\nu+\varepsilon}
 {(c_\ast+\varepsilon u+\varpi_\varepsilon)\sqrt{P_{r,\nu}(u)}}\,du
\]
is negative on $(r,1)$ and vanishes at both ends.  The two halves of the
oscillation give disjoint graphs of $z$ over the same monotone radial
coordinate, and the argument of Theorem~\ref{thm:rotating-torus} shows
that the meridian is a smooth embedded closed curve in $\{\rho>0\}$.  Its
surface of revolution is an embedded torus.  If $r\ge0$, $\mathcal
B_\varepsilon=0$ forces $\frac12x^2+\nu+\varepsilon$ to be negative somewhere
on $[r,1]$, and since it increases there,
\eqref{eq:split-torus-condition} holds.
\end{proof}

Above the parabola in Figure~\ref{fig:split-parameter-map}, the numerator of
$I_\varepsilon'$ is positive at both endpoints. Since $I_\varepsilon(r)=I_\varepsilon(1)=0$
for a closing pair, $I_\varepsilon$ changes sign from positive near $r$ to negative
near $1$. The two meridian branches thus cross at a positive radius, giving a
self-intersection of the swept surface.

Thus, in contrast with the boost family of
Theorem~\ref{thm:analytic-split-closure}, the elliptic branch has a time
return, and Proposition~\ref{thm:cylinder-torus-dichotomy} gives the
torus alternative, exactly as in the Heisenberg group.

\begin{figure}[!htbp]
 \centering
 \includegraphics[width=\textwidth]{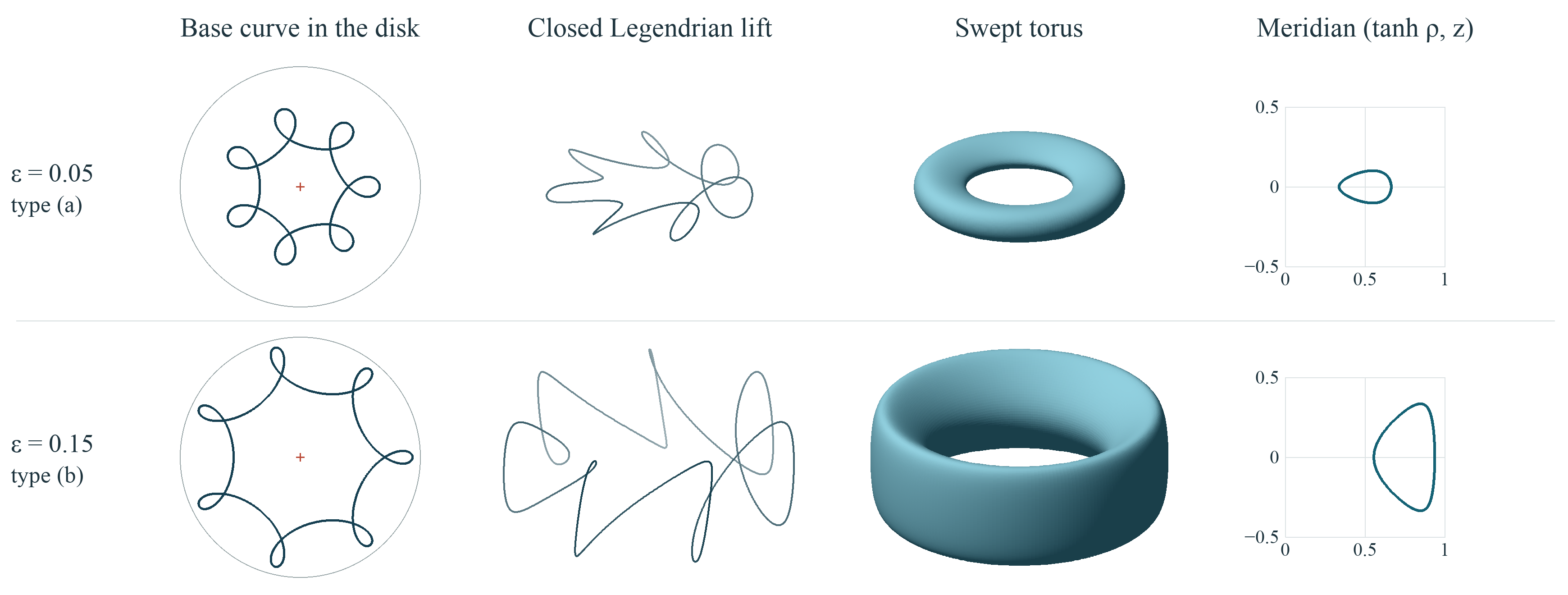}
 \caption{The two $6/7$ closing pairs of Figure~\ref{fig:split-parameter-map}.
 At $\varepsilon=0.05$ (top), $(r,\nu)\simeq(0.071987,-0.174047)$;
 at $\varepsilon=0.15$ (bottom), $(r,\nu)\simeq(-0.013181,-0.199912)$.
 The values of $\tfrac12 r^2+\nu+\varepsilon$ are approximately $-0.12146$
 and $-0.04983$, respectively, so Theorem~\ref{thm:split-torus} gives an
 embedded torus swept out by the lifted rotation in each case.
 Each row shows the base curve in the Poincar\'e disk centered at the
 rotation center (marked), its closed Legendrian lift, the swept torus,
 and its meridian. Both rows use the same scale in each column.
 The meridians are sections of the displayed tori in disk radius
 $\tanh\rho$ and fiber height $z$, with equal units on both axes
 and the same scale as the spatial views.}
 \label{fig:split-rotating-torus}
\end{figure}

\begin{remark}[Numerical closing pairs]
\label{rem:split-rotating-numerics}
Theorem~\ref{thm:split-rotating} gives existence for sufficiently small
$\varepsilon=4/b^2$. At the following sample values, numerical solution
of $\Theta/(2\pi)=p/n$ and $\mathcal B_\varepsilon=0$ gives the pairs
shown (rounded).
\[
\begin{array}{c|c|r|r}
p/n & \varepsilon & r & \nu\\ \hline
1/4  & 0.1920 & -0.1235988230 & -0.1996095384\\
1/2  & 0.2090 & -0.0855399767 & -0.2126695793\\
6/7  & 0.2363 & -0.0274120386 & -0.2367672876\\
9/10 & 0.2402 & -0.0193930469 & -0.2404883346
\end{array}
\]
For example, at $p/n=1/4$ the unrounded pair gives
$\varpi_\varepsilon^2\simeq4.033\times10^{-10}$ and
$c_*+\varepsilon r\simeq2.372\times10^{-5}$. Both quantities are
positive as required by Lemma~\ref{lem:split-rotating-data}.
The values of $\varepsilon$ in the table lie near the ends of the computed
elliptic branches: in every row $\varpi_\varepsilon^2$ is of order
$10^{-10}$, so these pairs show how far the branches extend rather than
typical closing pairs.  In the $1/4$ row
$\frac12r^2+\nu+\varepsilon\simeq2.9\times10^{-5}>0$, so
\eqref{eq:split-torus-condition} fails there.
The extent of the elliptic branches and their possible parabolic limits
remain open. Direct integration of \eqref{eq:hyperbolic-frame-matrix}
for $p/n=6/7$, $\varepsilon=0.05$ gives $\|\mathsf M^7-I\|<10^{-9}$
and absolute signed area below $10^{-8}$.
Figure~\ref{fig:split-rotating-torus} compares this torus with the one at $\varepsilon=0.15$.
  Unlike the
flat case (Proposition~\ref{prop:real-root-no-closure}), the branches pass
from type (a) to type (b) of Proposition~\ref{prop:general-cn-wave} at
$\varepsilon\approx0.10$--$0.15$; the lower row of
Figure~\ref{fig:split-rotating-torus} is of type (b).  Theorems
\ref{thm:split-rotating} and~\ref{thm:split-torus} do not distinguish the two
types.  In the split member the profile type of
Proposition~\ref{prop:general-cn-wave} and the type of the generator,
given by the sign of $\varpi_\varepsilon^2$, are independent, and there is no
analogue of Proposition~\ref{prop:real-root-no-closure}.  We treat only the
elliptic generators with $C\ne0$ and small $\varepsilon$, and the hyperbolic
generators with $C=0$ (Theorem~\ref{thm:analytic-split-closure}).
\end{remark}

\section{The swept surface and its pseudo-Hermitian invariants}
\label{sec:surface}

Let $f(s,t)=\gamma(s,t)$ be swept out by the Legendrian mKdV flow with $c_0=0$
in \eqref{eq:relative-equilibrium-velocity}, so that
$f_s=X$ and $f_t=\frac12k^2X+k_sJX+2k\xi$.  With $(X,JX,\xi)$ positively
oriented,
\begin{equation}
\label{eq:legendrian-normal}
 Z=f_s\times f_t=-2kJX+k_s\xi,\qquad N_H=\operatorname{pr}_HZ=-2kJX.
\end{equation}
Hence $k=k_s=0$ is rank degeneration, while on the regular set the
pseudo-Hermitian singular locus $N_H=0$ of Chiu--Huang--Lai~\cite{CHL2017}
is $\{k=0,\ k_s\ne0\}$.  Where $k\ne0$, $\eta(f_t)=2k\ne0$, so the
characteristic direction is $X$; with horizontal normal $n_H=JX$ the
pseudo-Hermitian mean curvature is
$H_p=\langle\nabla^{\TW}_XX,n_H\rangle=k$.  Thus $H_p$ itself satisfies
\eqref{eq:sasakian-mkdv}, and the $p$-area element is
$d\Sigma_p=|N_H|\,ds\,dt=2|H_p|\,ds\,dt$.  The densities $k$ and
$k^2$ are conserved, since
$(k^2)_t=(2kk_{ss}-k_s^2+\frac34k^4+2\lambda k^2)_s$; thus $\int2k\,ds$,
the signed horizontal Jacobian rate, and $\int H_p^2\,ds$ are constant.
If $k$ has fixed sign, the $p$-area rate $2\int|k|\,ds$ is constant as
well. It is also constant for any periodic traveling wave: translation
of $K(s-vt)$ preserves the integral of $|K|$ over a period, including
for the sign-changing figure-eight of Theorem~\ref{thm:closed-elliptic}.

\section{Dissipative and integrable Legendrian curve flows}
\label{sec:cim-comparison}

Pan--Sun~\cite{PanSun2021} study a parabolic, length-decreasing
curve-shortening flow of Legendrian curves in torsion-free
pseudo-Hermitian $3$-manifolds; the present flow is the Hamiltonian vector
field of the same length for $\Omega_\lambda$ and conserves it.  We use
\emph{dissipative} only in this sense and claim no gradient interpretation.
Calini--Ivey--Musso~\cite{CIM2024} develop the integrable hierarchy and its
spectral theory on the pseudo-Hermitian $3$-sphere; at $\lambda=2$,
Corollary~\ref{thm:sasaki-projection} recovers their compact mechanism.
The traveling-wave equation $K''+\frac12K^3+\alpha K=C$ is the
constrained-elastica equation after shifting the multiplier by the base
curvature, so the relative equilibria are the Euclidean, spherical, and
hyperbolic elastica of \cite{LangerSinger1984,Heller2014}.  The Legendrian
problem adds the Reeb holonomy condition and its consequences in
Remark~\ref{rem:relative-equilibrium-constant},
Proposition~\ref{thm:cylinder-torus-dichotomy}, and
Theorems~\ref{thm:split-rotating} and~\ref{thm:split-torus}.

\subsection*{Scope and outlook}
Does the strict increase in Remark~\ref{rem:flat-rotation-monotonicity}
extend to $-1<r<-1/2$? Then each prescribed angle in $(0,2\pi)$ would occur at exactly one
point of the flat zero-area graph. In the split member, open
questions concern closed relative equilibria
with hyperbolic generator and $C\ne0$, where the zero-area condition is no
longer forced by a reflection symmetry and one-parameter families are
expected, the range of continuation of the elliptic family in
Remark~\ref{rem:split-rotating-numerics}, and how the elliptic,
parabolic, and hyperbolic branches connect.  The positive $\dn$ branch
with $C=0$ is not treated here.  Its normalized constants are
$\nu=-(1+r^2)/4$ and $e_\ast=-r^2/8$, so
$\mathcal I=-4b^4\bigl((\nu+\varepsilon)^2-r^2/4\bigr)$ by
\eqref{eq:killing-casimir}, and its generator is elliptic, parabolic, or
hyperbolic according to the scale.  Its closed members, if any, would sweep
tori in the elliptic case and cylinders otherwise
(Proposition~\ref{thm:cylinder-torus-dichotomy}). Other questions include
a global hierarchy on Sobolev
completions, and
a spectral description of closed finite-gap relative equilibria.  We do not
study closed relative equilibria in the $3$-sphere beyond the general
framework of Section~\ref{sec:sasakian-family}; their closure is studied
in \cite{CIM2024}, and closed stationary solutions there belong to the
program announced in \cite{MussoUFRJ}.  A defocusing analogue would require a
Lorentzian horizontal metric, i.e.\ a para-Sasakian rather than a
pseudo-Hermitian setting.

\section*{Acknowledgements}

\noindent\textit{Funding.}
This work was supported in part by JSPS KAKENHI Grant Number JP22K03304.

\smallskip\noindent\textit{Competing interests.}
The author declares no competing interests.

\smallskip\noindent\textit{Use of generative AI.}
The author used the generative AI tools ChatGPT (OpenAI) and Claude
(Anthropic) in developing and checking proofs, in symbolic and numerical
calculations, and in drafting and revising the text. The author has
checked the statements, proofs, computations and references, and takes
full responsibility for the content of this paper.

\section*{Data availability statement}
The equations and parameter values used to compute the illustrated curves
are given in Section~\ref{sec:global-reconstruction} and in the figure
captions. Remark~\ref{rem:flat-rotation-monotonicity} reports the interval
bounds for the flat rotation number, and
Remark~\ref{rem:split-rotating-numerics} gives numerical closing pairs
for the hyperbolic model. The computations use no external datasets.

\end{document}